\documentclass[10pt]{article}
\usepackage{amssymb}
\usepackage{graphicx}
\def\qmod#1#2{{\hbox{}^{\displaystyle{#1}}}\!\big/\!\hbox{}_{
\displaystyle{#2}}}

\def\resto#1#2{{
#1\hskip 0.4ex\vline_{\hskip 0.4ex\raisebox{-0,5ex}
{{${\scriptstyle #2}$}}}}}
\def\at#1#2{
#1\hskip 0.25ex\vline_{\hskip 0.25ex\raisebox{-1.5ex}
{{$\scriptstyle#2$}}}}

\newfam\msbfam

\font\tenmsb=msbm10

\font\sevenmsb=msbm10 at 7pt
\font\fivemsb=msbm10 at 5pt

\textfont\msbfam=\tenmsb
\scriptfont\msbfam=\sevenmsb
\scriptscriptfont\msbfam=\fivemsb

\def\C{{\mathbb C}}

\def\F{{\mathbb F}}

\def\H{{\mathbb H}}

\def\N{{\mathbb N}}

\def\P{{\mathbb P}}
\def\Q{{\mathbb Q}}
\def\R{{\mathbb R}}
\def\S{{\mathbb S}}

\def\V{{\mathbb V}}
\def\Z{{\mathbb Z}}

\def\union{\mathop{\bigcup}}
\def\qed {\hfill\vrule height6pt width6pt depth0pt \bigskip}

\def\map{\longrightarrow}
\def\textmap#1{\mathop{\vbox{\ialign{
                                 ##\crcr
     ${\scriptstyle\hfil\;\;#1\;\;\hfil}$\crcr
     \noalign{\kern 1pt\nointerlineskip}
     \rightarrowfill\crcr}}\;}}

\def\textlmap#1{\mathop{\vbox{\ialign{
                                 ##\crcr
     ${\scriptstyle\hfil\;\;#1\;\;\hfil}$\crcr
     \noalign{\kern-1pt\nointerlineskip}
     \leftarrowfill\crcr}}\;}}

\newfam\meuffam

\font\tenmeuf=eufm10
\font\sevenmeuf=eufm7
\font\fivemeuf=eufm5

\textfont\meuffam=\tenmeuf
\scriptfont\meuffam=\sevenmeuf
\scriptscriptfont\meuffam=\fivemeuf

\def\germ{\fam\meuffam\tenmeuf}

\def\g{{\germ g}}
\def\hg{{\germ h}}

\newtheorem{sz}{Satz}[section]
\newtheorem{thry}[sz]{Theorem}
\newtheorem{pr}[sz]{Proposition}
\newtheorem{re}[sz]{Remark}
\newtheorem{co}[sz]{Corollary}
\newtheorem{dt}[sz]{Definition}
\newtheorem{lm}[sz]{Lemma}

\begin{document}
\def\Pr{{\rm Pr}}
\def\tr{{\rm Tr}}
\def\End{{\rm End}}
\def\Aut{{\rm Aut}}
\def\Spin{{\rm Spin}}
\def\U{{\rm U}}
\def\SU{{\rm SU}}
\def\SO{{\rm SO}}
\def\PU{{\rm PU}}
\def\GL{{\rm GL}}
\def\spin{{\rm spin}}
\def\su{{\rm su}}
\def\so{{\rm so}}
\def\ub{\underbar}
\def\pu{{\rm pu}}
\def\Pic{{\rm Pic}}
\def\Iso{{\rm Iso}}
\def\NS{{\rm NS}}
\def\deg{{\rm deg}}
\def\Hom{{\rm Hom}}
\def\Aut{{\rm Aut}}
\def\h{{\germ h}}
\def\Herm{{\rm Herm}}
\def\Vol{{\rm Vol}}
\def\pf{{\bf Proof: }}
\def\id{{\rm id}}
\def\i{{\germ i}}
\def\im{{\rm im}}
\def\rk{{\rm rk}}
\def\ad{{\rm ad}}
\def\h{{\bf H}}
\def\coker{{\rm coker}}
\def\dbar{\bar{\partial}}
\def\Lo{{\Lambda_g}}
\def\niq{=\kern-.18cm /\kern.08cm}
\def\Ad{{\rm Ad}}
\def\RSU{\R SU}
\def\ad{{\rm ad}}
\def\dva{\bar\partial_A}
\def\da{\partial_A}
\def\p{\partial\bar\partial}
\def\sp{\Sigma^{+}}
\def\sm{\Sigma^{-}}
\def\spm{\Sigma^{\pm}}
\def\smp{\Sigma^{\mp}}
\def\Tors{{\rm Tors}}
\def\st{{\rm st}}
\def\s{{\rm s}}
\def\oo{{\scriptstyle{\cal O}}}
\def\ooo{{\scriptscriptstyle{\cal O}}}
\def\sw{Seiberg-Witten }
\def\pa{\partial_A\bar\partial_A}
\def\Dr{{\raisebox{0.15ex}{$\not$}}{\hskip -1pt {D}}}
\def\gr{{\scriptscriptstyle|}\hskip -4pt{\g}}
\def\subsetint{{\  {\subset}\hskip -2.45mm{\raisebox{.28ex}
{$\scriptscriptstyle\subset$}}\ }}
\def\ra{\rightarrow}
\def\pst{{\rm pst}}
\def\sst{{\rm sst}}

\def\kod{{\rm kod}}
\def\degmax{{\rm degmax}}
\def\red{{\rm red}}
\def\ASD{{\rm ASD}}

\title{Harmonic sections in sphere bundles, normal neighborhoods of reduction loci, and  instanton moduli spaces on definite 4-manifolds}
\author{Andrei Teleman}
\maketitle

\begin{abstract}
In the first part of the paper we prove an existence theorem for gauge invariant $L^2$-normal neighborhoods of the reduction loci in the space ${\cal A}_a(E)$ of oriented connections on a fixed  Hermitian 2-bundle $E$.  We use this to obtain results on the topology of the moduli space ${\cal B}_a(E)$  of (non-necessarily irreducible) oriented connections, and to study the Donaldson $\mu$-classes globally  around the reduction loci. In this part of the article we use essentially  the concept of harmonic section in a sphere bundle with respect to an Euclidean connection.

Second, we concentrate on moduli spaces of instantons on definite 4-manifolds with arbitrary  first Betti number. We prove strong generic regularity  results which imply  (for bundles with ``odd" first Chern class) the existence of a  {\it connected}, dense  open  set of ``good" metrics for which all the reductions in the Uhlenbeck compactification of the moduli space are simultaneously  regular. These results can be used to define new Donaldson type invariants for definite 4-manifolds. The  idea behind this construction is to notice that, for a good metric $g$, the geometry of the instanton moduli spaces     around the reduction loci is always the same, independently of the choice of $g$.  The connectedness of the space of good metrics is important, in order to prove that no wall-crossing phenomena (jumps of invariants) occur.  Moreover, we notice that, for low instanton numbers, the corresponding moduli spaces  are a priori compact and contain no reductions at all so, in these cases, the existence of well-defined Donaldson type invariants is obvious. Note that, on the other hand, there seems to be no way to introduce well defined  numerical Seiberg-Witten invariants for definite 4-manifolds: the construction proposed in \cite{OT2} gives a $\Z$-valued function defined on a countable set of chambers.

The natural question is to decide whether these new Donaldson type invariants    yield essentially new {\it  differential} topological information on the base manifold have, or have a purely topological nature.  
\end{abstract}

\tableofcontents

\setcounter{section}{-1}
\section{Introduction}

The main goal of this article is to study moduli spaces of instantons over 4-manifolds with negative definite intersection form. The vanishing of $b_+$ has an important consequence on the geometry of the instanton moduli spaces: all line bundles admit ASD connections (with respect to any metric) hence, as soon as rank 2-bundle $E$  splits topologically, the corresponding instanton moduli space  will always  contain reductions. In other words, ``one cannot get rid of reductions by perturbing the metric".  On the other hand, our main applications will concern 4-manifolds with $b_1>0$, and for such manifolds the spaces of reductions are positive dimensional.  Therefore, it is very important to study carefully the   {\it global} geometry of the moduli space of connections around the loci of reductions. This will be our first goal.

To be more precise, let $E$ be rank 2-Hermitian bundle on a 4-manifold $M$, and denote $D:=\det(E)$, $d=c_1(D)$. Consider the affine space ${\cal A}_a(E)$  of connections $A$ on   $E$ which induce a fixed connection $a$ on  $D$, and the moduli space ${\cal B}_a(E)={\cal A}_a(E)/{\cal G}_E$, where ${\cal G}_E$ is the gauge group $\Gamma(SU(E))$.

Let $l\in  H^2(X,\Z)$ such that $l(d-l)=c_2(E)$ and consider the set $\lambda=\{l,d-l\}$ (which has either one or two elements).  We denote by ${\cal A}_a^\lambda(E)$ the subspace of connections  $A\in {\cal A}_a(E)$, which  are simply reducible of type $\lambda$, i.e. which admit only two parallel line subbundles  whose Chern classes are $l$, $d-l$. Such a connection will be called $\lambda$-reducible.  ${\cal A}_a^\lambda(E)$ becomes a (locally closed) submanifold of the affine space  ${\cal A}_a(E)$ (after suitable Sobolev completions).  Our first problem is the construction of a gauge invariant $L^2$-normal neighborhood of this submanifold. More precisely, we will show that, denoting by $N^\lambda$ the $L^2$-normal bundle of ${\cal A}_a^\lambda(E)$, the restriction  of the natural map
$\nu:N^\lambda\to {\cal A}_a(E)$ to a sufficiently small gauge invariant neighborhood ${\cal U}^\lambda$ of the zero section  is a diffeomorphism on its image. Moreover, the neighborhood ${\cal U^\lambda}$ is defined by a inequality of the form $\|\alpha\|_{L^\infty}\leq\varepsilon(A)$ (on the fibre $N_A$), where  the assignment $A\mapsto \varepsilon(A)$ is gauge invariant and  continuous (with respect to a sufficiently fine Sobolev topology on the space of connections).

Although this statement is very natural, the proof is not easy. The difficulty comes from the fact that infinite dimensional manifolds are not locally compact. Even the fact the $\nu$ is injective on a neighborhood of the zero section is not trivial. The main difficulty is to characterize in a convenient way the connections which are ``close" to the reducible locus  ${\cal A}_a^\lambda(E)$, i.e. which are ``almost" $\lambda$-reducible.

 Our argument is based on the following idea: A reduction  $A\in {\cal A}_a^\lambda(E)$ admits a parallel section  in the sphere bundle $S(su(E))$. A connection which is close to being reducible should admit an {\it energy minimizing harmonic section} in this sphere bundle. The precise meanings of  the words ``energy"  and ``harmonic" are the following: we associate to any connection $A\in {\cal A}_a(E)$ the energy functional
 $$E_A(u):=\|d_A(u)\|_{L^2}^2\ ,\ E_A:\Gamma(S(su(E)))\map \R_{\geq 0}\ .
 $$
 and we agree to call the critical points of this functional $A$-harmonic sections. After proving these results about normal neighborhoods of reduction loci, we realized that our problem can be naturally  generalized in the following way: For an Euclidean bundle $F$ on an arbitrary compact Riemannian manifold, construct a gauge equivariant $L^2$-normal neighborhood of the space of Euclidean connections $d_A:A^0(F)\to A^1(F)$ with  1-dimensional kernel. Therefore, in  the first section we will treat this more general problem, which is of independent interest; the results concerning reducible Hermitian connections will be easily deduced as applications, taking $F=su(E)$.  
 
 The second section deals with the first consequences of our existence results for gauge  invariant $L^2$-normal neighborhoods: a global description of the universal $SO(3)$ bundle, explicit formulae for the Donaldson  $\mu$-classes around the reduction loci and an explicit description of the homotopy type of the whole moduli space ${\cal B}_a(E)$. This space is much more complex than its open subspace ${\cal B}^*_a(E)$, which plays a central role in classical gauge theory and whose rational homotopy type has been described in \cite{DK}. Our  description of  ${\cal B}_a(E)$ yields an easy method to compute the cohomology of this space using the Mayer-Vietoris exact sequence.
 
 The third section is dedicated to the geometry of the instanton  moduli space around a reduction locus.  Using our existence results for normal neighborhoods, we show that {\it in a neighborhood of the reduction  locus associated with a fixed topological splitting of the bundle, the instanton moduli problem reduces to an abelian moduli problem}, which is very much similar to the Seiberg-Witten one. These abelian equations read
 \begin{equation} 
\left\{
\begin{array}{ccc}
F_b^+-\frac{1}{2}F_a^+&=&(\alpha\wedge\bar\alpha)^+\ .\vspace{2mm}  \\  
(d_{b^{\otimes 2}\otimes a^\vee}^*,d_{b^{\otimes 2}\otimes a^\vee}^+)\alpha&=&0\ ,
\end{array}
\right.
\end{equation}
which are equations for a pair $(b,\alpha)$, where $b$ is a connection  on fixed Hermitian line bundle $L$ and $\alpha\in A^1(L^{\otimes 2}\otimes D^\vee)$.  This result provides  simple descriptions  of the linear spaces of  harmonic spaces of the deformations elliptic complexes at the reductions. It is very important to have a {\it global} description of this linear space. Similar results are stated for the loci of {\it twisted reductions} (i.e. instantons  which are locally reducible, but globally irreducible).

Our most important technical results are obtained in the last subsection of section 3: we prove strong generic regularity theorems at the reductions. Regularity at the reductions is an old, classical problem in gauge theory (see \cite{FU}, \cite{DK}) and one might wonder whether there are  still unsolved questions on this problem. The point is that Freed-Uhlenbeck's generic regularity result is not sufficient for our purposes. For our purposes, we need a  {\it connected}, dense open  set of good metrics (metrics for which all reductions in the Uhlenbeck compactification are regular). 
Our proof has two steps:
\\ 
\\
1. Define a  connected, dense open connected open of {\it admissible} ${\cal C}^r$-metrics, for which the vanishing loci  of the harmonic representatives  of the classes 
$$2l-d\ ,\ l\cdot (d-l)\leq c_2(E)$$
  have good geometric properties (see section \ref{admissmet} in the Appendix). Our admissibility condition is very natural: we require that the rank of the intrinsic derivative of  these harmonic representatives at any vanishing point is    at least 2. The difficulty is to check that this condition defines indeed a  connected, dense set of metrics.

In particular, for our  admissible metrics, the  vanishing loci of all  these harmonic representatives  have Hausdorff dimension at most 2. Note that this holds  for {\it any} ${\cal C}^\infty$-metric by a result of B\"ar  \cite{B}, but this result does not appear to generalize for ${\cal C}^r$-metrics. 
For a metric for which this Hausdorff dimension  bound holds, the statement in   Lemma 4.16   \cite{FU} holds, making possible the second step. This stronger version of Lemma 4.16   \cite{FU} is proved in detail in the Appendix (Corollary \ref{fub}).
\\
 \\ 
2. Regard the linear space formed by the {\it second} harmonic spaces at the reductions  as the moduli space associated with an {\it abelian} moduli problem, and prove a transversality theorem for this moduli problem with respect to variations of $g$ (in the space of admissible metrics).  In other words, we will prove that the parameterized moduli space (obtained by letting the metric vary in the set of admissible metrics) is smooth away of the zero-section.
Next we  show that  -- for a bundle $E$ with odd Chern class --    the projection map from the $\C^*$-quotient of this parameterized moduli space (minus the zero-section) on the space of admissible metrics  is Fredholm of {\it negative} index $\leq -2$, hence its image has  connected complement.  \\

Section 4 deals with applications of our results. First we prove a simple geometric property of a particular  instanton moduli spaces on a 4-manifold which has the homology type of a class VII surface with  $b_2=2$: the two circles  of reductions belong always (for any metric!) to the same connected component. We continue with the construction of the new class of Donaldson invariants.

\section{Harmonic sections, parallel sections} 

\subsection{Harmonic sections in sphere bundles}

Let $(M,g)$ be a compact oriented Riemannian $n$-manifold and  let $F$ be a real  rank $r$ vector bundle  on $M$  endowed with an inner product, and denote by $S(F)$ the unit sphere bundle of $F$.   Let $A$ be an Euclidean connection on $F$.
The energy functional on the space of of sections $\Gamma(S(F))$ is defined by
$$E_A(u)=\| d_A (u)\|^2_{L^2}=\langle d_A u,d_A u\rangle^2\ .
$$
The critical points of this functional will be called $A$-harmonic sections. 
\begin{pr} A section $u\in \Gamma(S(F))$ is harmonic if and only if it satisfies the Euler-Lagrange equation
$$d_A^*d_A u-|d_A u|^2 u=0\ .
$$
\end{pr}
\pf The section  $u\in \Gamma(S(F))$ is a critical point  of $E_A$ if and only if 
$$\langle d_A^*d_A u, v\rangle_{L^2}=0
$$
for every $v\in T_u(\Gamma(S(F)))$.  This happens if and only if  there exists a real function $\varphi$ such that $d_A^*d_A u=\varphi u$. On the other hand, the well-known identity
$$\frac{1}{2}\Delta |u|^2=(d_A^*d_A u,u)- |d_A u|^2
$$
shows that one must have $\varphi=|d_A u|^2$.
\qed

For a fixed connection $A$ the theory of $A$-harmonic sections is very much similar to the theory of sphere valued harmonic maps. In particular, one has a parabolic evolution equation given by the gradient flow of the functional $E_A$, and using this equation and its long-time convergence properties, one can study the existence of a $A$-harmonic representative in a given homotopy class of sections.\\

Two vectors $a$, $b$ of an Euclidean vector space $V$ define an endomorphism $a\wedge_V b\in so(V)$ given by $a\wedge_V b(h)=(a,h)b-(b,h)a$.  The subspace
$$\{a\wedge_V b|\  b\in V\}\subset so(V)
$$
coincides with the  orthogonal complement  $so(a^\bot)^\bot$ of $so(a^\bot)$ in $so(V)$, and the map $a^\bot \ni b\mapsto a\wedge_V b\in so(a^\bot)^\bot$ is an isomorphism.
Similarly, for a section $u\in A^0(F)$ and a form $v\in A^1(F)$ we obtain a form $u\wedge_F v\in A^1(so(F))$.
\begin{re}\label{A0} Let $A \in{\cal A}(F)$, $u\in\Gamma(S(F))$ and $A_0:=A-u\wedge_F d_A u$. Then
\begin{enumerate}
\item The section  $u$ is $A_0$-parallel.
\item  $A_0$ is the unique connection making $u$ parallel with $A-A_0\in A^1(so(u^\bot)^\bot)$.
\end{enumerate}
\end{re}
\pf  1. One has 
$$d_{A_0} u=d_A u-(u\wedge_F d_A u)(u)= 
d_A u-(u,u)d_A u+(u,d_Au)u=0\ ,
$$
because, since $(u,u)\equiv 1$, one has $(u,d_Au)=0$. \\

2. If $A-A_0\in A^1(so(u^\bot)^\bot)$, there exists a  unique form $v\in A^1(u^\bot)$ such that $A-A_0=u\wedge_F v$. The relation $d_{A_0} u=0$ is equivalent to $d_A u=v$ .
\qed
\begin{pr} Let $u\in \Gamma(S(F))$ and $A\in{\cal A}(F)$. The following conditions are equivalent
\begin{enumerate}
\item The section $u$ is $A$-harmonic.
\item Putting  $A_0:=A-u\wedge_F d_A u$, one has $d_{A_0}^* (A-A_0)=0$.
\end{enumerate}
\end{pr}
\pf  We write locally $d_A u=\nabla^A_{e_i}u e^i$ with respect to a local orthonormal frame  $(e_i)$. Therefore, taking into account that $u$ is $A_0$ parallel  and   $(u,\nabla_{e_i}^Au)=0$, one has 
$$d_{A_0}^* (A-A_0)=(d_{A_0})^* (u\wedge_F d_A u)=u\wedge_F (d_{A_0})^*d_A u=
u\wedge_F (d_{A}-u\wedge_F d_A u)^*d_A u=
$$
$$=u\wedge_F[ d_A^* d_A u+*(u\wedge_F d_A u)*d_Au]=u\wedge_F[ d_A^* d_A u+*(u\wedge_F \nabla^A_i u e^i)*\nabla^A_j u e^j]
$$
$$=u\wedge_F[ d_A^* d_A u- (\nabla^A_i u, \nabla^A_j u) u *(e^i\wedge * e^j)]=u\wedge_F[ d_A^* d_A u- |d_Au|^2 u]\ .$$
Since $d_A^* d_A u- |d_Au|^2 u\in A^0(u^\bot)$, the vanishing of $u\wedge_F[ d_A^* d_A u- |d_Au|^2 u]$ is equivalent to the vanishing of $d_A^* d_A u- |d_Au|^2 u$.
\qed

 The statement of the proposition above can be formulated as follows
 \begin{co} If the connection $A$ admits a harmonic section $u\in \Gamma(S(F))$, then it is in Coulomb gauge with respect to a connection $A_0$ for which $u$ is a parallel section. 
 \end{co}
 
 \subsection{A normal neighborhood of the locus of Euclidean connections with 1-dimensional kernel}\label{nn}
 
 We define the locally closed subspace
$${\cal A}'(F):=\{A\in{\cal A}(F)|\ \dim(\ker(d_A))=1\}\subset {\cal A}(F)\ .
$$
Let ${\cal S}(F)$ be the space of trivial rang 1 subbundles of $F$. This space can be naturally identified with $\Gamma(S(F))/\{\pm 1\}$, and becomes a smooth manifold after suitable Sobolev completions. Two rang 1 subbundles which are sufficiently ${\cal C}^0$-close have isomorphic complements, so they are conjugate modulo the action of the gauge group $\Aut(F)=\Gamma(SO(F))$. In other words, the action of $\Aut(F)$ on ${\cal S}(F)$  is locally transitive. One has an obvious surjective map 
$$w:{\cal A}'(F)\to {\cal S}(F)\ ,\ A\mapsto \langle\ker(d_A)\rangle\ ,$$
where $\langle\ker(d_A)\rangle=\im({\rm ev}:\ker(d_A)\times M\to E)$ stands for the subbundle generated by the line $\ker(d_A)$.     We  will use the  subscript $(\cdot )_k$ (for $k\in\N$ sufficiently large) to denote Sobolev completion with respect to the $L^2_k$-norm.
\begin{pr}\label{subm}  The subset  ${\cal A}'(F)_k\subset{\cal A}(F)_k$ is a submanifold, and the surjection  $w:{\cal A}'(F)_k\to {\cal S}(F)_{k+1}$ is a submersion.  
\end{pr}
\pf We omit Sobolev indices to save on notations. Let $A\in {\cal A}'(F)$ and $\lambda=\langle \ker(d_A)\rangle\in{\cal S}(F)$.       The fibre $w^{-1}(\lambda)$ is obviously  an affine subspace of ${\cal A}(F)$ which can be identified with ${\cal A}(u^\bot)$.  The stabilizer $H$ of  $\lambda$ is a closed Lie subgroup $H$ of the gauge group $\Aut(F)$ whose Lie algebra $\hg$ can be identified  with $A^0(so(u^\bot))$, so it has a   topological complement $\hg^\bot=A^0(so(u^\bot)^\bot)$. The restriction of the map $\hg^\bot\times w^{-1}(\lambda)\to {\cal A}(F)$ given by $(h,B)\mapsto \exp(h)(B)$ to a sufficiently small open neighborhood of $(0,A)$ is an embedding which parameterizes  a neighborhood  of $A$ in ${\cal A}'(F)$.  This gives the submanifold structure of ${\cal A}'(F)$.
Using the local transitivity of the gauge action on  ${\cal S}(F)$, it follows that $h\mapsto \exp(h)(A)$  defines a local slice of $w$ at $A$, proving that $w$ is a submersion at $A$.\qed

The main goal of this section is to prove that 
\begin{thry}
The submanifold ${\cal A}'(F)\subset {\cal A}(F)$ has a tubular, gauge invariant, $L^2$-normal neighborhood.
\end{thry}
The idea is very simple and natural: construct a neighborhood of ${\cal A}'(F)$ consisting of connections which admit an (up to sign) unique harmonic, energy minimizing section.\\

Let $A_0\in {\cal A}'(F)$, and let $u$ one of the two generators in $\Gamma(S(F))$ of the line $\ker(d_{A_0})$.  The splitting $F=\langle u \rangle\oplus u^\bot$ defines an $A_0$-parallel splitting
$$so(F)=so(u^\bot)\oplus so(u^\bot)^\bot=so(u^\bot)\oplus u\wedge_Fu^\bot\ , 
$$
which gives an $L^2$-orthogonal decomposition
$$A^1(so(F))=A^1(so(u^\bot))\oplus (d_{A_0})[A^0(so(u^\bot)^\bot)]\oplus $$
$$\oplus \ker[(d_{A_0})^*:A^1((so(u^\bot)^\bot)\to A^0(so(u^\bot)^\bot)]\ .
$$

The geometric interpretations of the three factors in the decomposition above are the following: the space $A^1(so(u^\bot))$ is the tangent space of the fibre $w^{-1}(\langle u\rangle)$ at $A$, $n_A:=(d_{A_0})[A^0(so(u^\bot)^\bot)]$ is the normal space at $A_0$ of this fibre in the submanifold ${\cal A}'(F)$, whereas the space 
$$N_{A_0}:=\ker\left[(d_{A_0})^*:A^1((so(u^\bot)^\bot)\to A^0(so(u^\bot)^\bot)\right]$$
 is the normal space at $A_0$  of  the submanifold  ${\cal A}'(F)$ in the space of connections ${\cal A}(F)$.
 
Let  $N\to {\cal A}'(F)$ be the normal vector bundle of the submanifold ${\cal A}'(F)$, whose  fibre $A_0\in {\cal A}'(F)$ is just the space  $N_{A_0}$ defined above. One has a natural map $\nu:N\map {\cal A}(F)$ given by
$$\nu(A_0,\alpha)={A_0}+\alpha\ ,\ \forall A_0\in{\cal A}'(F)\ ,\ \alpha\in N_{A_0}\ ,
$$
which is obviously a local isomorphism at every point of the form $(A_0,0)$. The map $\nu$ is equivariant with respect to the natural gauge actions on $N$ and ${\cal A}(F)$. For $\zeta\in A^0(so(F))={\rm Lie}(\Aut(F))$ we denote by $\zeta^\#$ the tangent field (the infinitesimal transformation) of $N$ associated with $\zeta$ and by $\zeta_\#$ the corresponding tangent field of ${\cal A}(F)$. The equivariance property of $\nu$ implies
\begin{equation}\label{equivformula}
d(\nu)(\zeta^\#_{(A_0,\alpha_0)})=[\zeta_\#]_{A_0+\alpha_0}=d_{A_0+\alpha_0}(\zeta)=d_{A_0}(\zeta)+[\alpha_0,\zeta]\ .
\end{equation}

\begin{lm}\label{open} Let $A_0\in{\cal A}'(F)$. There exists $\varepsilon(A_0)>0$ such that the differential $d_{(A_0,\alpha_0)}\nu$ at  $(A_0,\alpha_0)$ is an isomorphism for every $\alpha_0\in N_{A_0}$ with $\|\alpha_0\|_{L^\infty}\leq\varepsilon(A_0)$. The assignment $A_0\mapsto \varepsilon(A_0)$ can be chosen to continuous and gauge invariant.
\end{lm}
\pf The submersion $w:{\cal A}'(F)\to{\cal S}(F)$ induces a submersion $\tilde w:N\to {\cal S}(F)$.  Fix $(A_0,\alpha_0)\in N$ and choose   $u_0\in\ker (d_{A_0})\cap\Gamma(S(F))$.  The subspace
$$\{\zeta^\#_{(A_0,\alpha_0)}|\ \zeta\in A^0(so(u_0^\bot)^\bot)\}
$$
is a $\tilde w$-horizontal space at  $(A_0,\alpha_0)$, i.e. a topological complement of the vertical tangent space $T_{(A_0,\alpha_0)}(\tilde w^{-1}(\langle u_0\rangle)$. This complement   is isomorphic with $A^0(so(u_0^\bot)^\bot)$. A vertical   tangent vector 
$$v'\in T_{(A_0,\alpha_0)}(\tilde w^{-1}(\langle u_0\rangle)$$
 can be written as $v'=(b,\beta)$, where $b\in A^1(so(u_0^\bot))$ and $\beta\in A^1(so(u_0^\bot)^\bot)$ satisfy the equation 
 $$d_{A_0}^*\beta-*[b\wedge*\alpha_0]=0\ .
 $$
(obtained by derivating the relation $d_{A_0}^*\alpha_0=0$ in the direction $(\dot A_0,\dot{\alpha_0})=(b,\beta)$). Therefore, one has an isomorphism
$$T_{(A_0,\alpha_0)}(N)= \left\{(b,\beta)\in A^1(so(u_0^\bot))\oplus A^1(so(u_0^\bot)^\bot)\vline\  d_{A_0}^*\beta-*[b\wedge*\alpha_0]=0\right\}\oplus $$
$$\oplus A^0(so(u_0^\bot)^\bot)\ .
$$
Using (\ref{equivformula}) one obtains
$$d\nu(b,\beta,\zeta)=b+\beta + d_{A_0}(\zeta)+[\alpha_0,\zeta]\ .
$$
The statement follows now directly from Lemma \ref{isomop} below.
\qed
\begin{lm}\label{isomop} If $\|\alpha_0\|_{L^\infty}$ is sufficiently small, then the operator
$$A^1(so(u_0^\bot))_k\oplus A^1(so(u_0^\bot)^\bot)_k\oplus A^0(so(u_0^\bot)^\bot)_{k+1}\stackrel{P}{\to}   A^1(so(F))_k\oplus A^0(so(u_0^\bot)^\bot)_{k-1}$$
given by
$$ (b,\beta,\zeta)\mapsto   (b+\beta + d_{A_0}(\zeta)+[\alpha_0,\zeta],d_{A_0}^*\beta-*[b\wedge*\alpha_0]  $$
 is an isomorphism.
\end{lm}
\pf We omit as usually the Sobolev indexes. Using the decomposition
$$A^1(so(F))=A^1(so(u_0^\bot))\oplus A^1(so(u_0^\bot)^\bot)\ ,
$$
the operator $P$ can be written as
$$P=\left[
\begin{array}{ccc}
\id&0&[\alpha_0,\cdot]\\
0&\id&d_{A_0}\\  
-*[\cdot\wedge*\alpha_0]&d_{A_0}^*&0
\end{array}
\right]
$$
If $(b,\beta,\zeta)\in\ker (P)$, one obtains easily 
$$b=-[\alpha_0,\zeta]\ ,\ d_{A_0}^*\beta=-d_{A_0}^*d_{A_0}\zeta=*(b\wedge*\alpha_0)=*(-[\alpha_0,\zeta]\wedge*\alpha_0)$$
\begin{equation}\label{eq}
\Delta_{A_0}\zeta =*[[\alpha_0,\zeta]\wedge*\alpha_0]
\end{equation}

On the other hand, since $u_0$ is parallel and $\dim(\ker(d_{A_0}))=1$,  one has  
$$\ker(d_{A_0}: A^0(su(u_0^\bot)^\bot)\to A^1(su(u_0^\bot)^\bot))=$$
$$=u_0\wedge_F\ker(d_{A_0}:A^0(u_0^\bot)\to A^1(u_0^\bot))=\{0\}
$$
so the minimal eigenvalue $\eta(A_0)$ of $\Delta_{A_0}: A^0(su(u_0^\bot)^\bot)\to A^0(su(u_0^\bot)^\bot)$ is positive. Therefore (\ref{eq}) implies $\zeta=0$ (hence $\ker(P)=0$) as soon as $\|\alpha_0\|^2_{L^\infty}$ is sufficiently small with respect to $\frac{1}{\eta(A_0)}$. For surjectivity, note that the equation  
\begin{equation}\label{surj}
P(b,\beta,\zeta)=(c,\gamma,\eta)
\end{equation}
 becomes
$$b+[\alpha_0,\zeta]=c\ ,\ \beta+d_{A_0}\zeta=\gamma\ ,\ d_{A_0}^*\beta-*[b\wedge*\alpha_0]=\eta\ .
$$
Consider first the weaker equation (for the single unknown $\zeta$)
$$\Delta_{A_0}\zeta-*[[\alpha_0,\zeta]\wedge*\alpha_0]=-*[c\wedge*\alpha_0]-\eta+d_{A_0}^*\gamma\ ,
$$
which is uniquely solvable, as soon as  $\|\alpha_0\|^2_{L^\infty}$ is sufficiently small, so one gets a unique solution $\zeta\in A^0(so(u_0^\bot)^\bot)_{k+1}$ for a triple $(c,\gamma,\eta)$ of Sobolev type $(L^p_k,L^p_k,L^p_{k-1})$. Then set $\beta:=\gamma-d_{A_0}\zeta$, $b:=c-[\alpha_0,\zeta]$, and we get a solution of the equation (\ref{surj}) of Sobolev type $(L^p_k,L^p_k,L^p_{k+1})$.
\qed
\begin{lm}\label{absmin} Let $A_0\in{\cal A}'(F)$ and $u_0\in\ker (d_A)\cap\Gamma(S(F))$. There exists $\varepsilon(A_0)>0$ such that for every $\alpha\in A^1(so(u^\bot_0)^\bot)$ with $d_{A_0}^*(\alpha)=0$ and $\|\alpha\|_{L^\infty}\leq\varepsilon(A_0)$ the energy functional $E_{A_0+\alpha}$ on the space $\Gamma(S(F))$ obtains its absolute minimum  at   $\pm u_0$ and only at these sections. Moreover, the assignment ${\cal A}'(F)\ni A_0\mapsto\varepsilon(A_0)$ can be chosen to be continuous and gauge invariant.
\end{lm}
\pf Consider a section $u\in \Gamma(S(F)))$,   put $v:=u-u_0$, $v':= u+u_0$, and set $\alpha:=u_0\wedge_F a$, where $a\in A^1(u_0^\bot)$. This implies $|\alpha|^2=2|a|^2$. One has
$$\|d_{A_0+\alpha} u\|^2-\|d_{A_0+\alpha} u_0\|^2=\|d_{A_0+\alpha} v\|^2+2\langle d_{A_0+\alpha} u_0,d_{A_0+\alpha} v\rangle=
$$
$$
=\|d_{A_0}v\|^2+\|\alpha(v) \|^2+2\langle d_{A_0}v,\alpha(v)\rangle+2\langle \alpha(u_0),d_{A_0}v\rangle+2\langle \alpha(u_0),\alpha(v)\rangle\ .
$$
But
$$d_{A_0}^*(\alpha(u_0)=-*d_{A_0}(*\alpha(u_0))=-*d_{A_0}(*\alpha)(u_0)\pm*(*\alpha\wedge d_{A_0}u_0)=0$$
(because $d_{A_0}^*(\alpha)=0$ and $d_{A_0} u_0=0$). On the other hand it holds pointwise
$$(\alpha(u_0),\alpha(v))=((u_0\wedge_F a)(u_0),(u_0\wedge_F a)(v))=$$
$$
=(a,(u_0,v)a-(a,v)u_0)=|a|^2(u_0,v)=\frac{1}{2}|\alpha|^2(u_0,v),
$$
$$1=|u_0+v|^2=|u_0|^2\ ,\  2(u_0,v)=-|v|^2\ .
$$
Therefore
$$\|d_{A_0+\alpha} u\|^2-\|d_{A_0+\alpha} u_0\|^2=\|d_{A_0}v\|^2+\|\alpha(v)\|^2+2\langle d_{A_0}v,\alpha(v)\rangle-\frac{1}{2}\int_M|\alpha|^2|v|^2$$
$$\geq \|d_{A_0}v\|^2-c_1 \sup|\alpha| \|d_{A_0}v\|_{L^2}\|v\|_{L^2}-c_2\sup|\alpha|^2\|v\|_{L^2}^2
$$
Since the same computation also applies to $-u_0$, we get
\begin{equation}\label{ineq1}
E_{A_0+\alpha}(u)-E_{A_0+\alpha}(u_0)\geq \|d_{A_0}v'\|^2-c_1 \sup|\alpha| \|d_{A_0}v'\|\|v'\|-c_2\sup|\alpha|^2\|v'\|^2
$$
$$
E_{A_0+\alpha}(u)-E_{A_0+\alpha}(u_0)\geq \|d_{A_0}v\|^2-c_1 \sup|\alpha| \|d_{A_0}v\|\|v\|-c_2\sup|\alpha|^2\|v\|^2
\end{equation}

Since $\|u\|_{L^2}^2=\|u_0\|_{L^2}^2=\Vol(M)$, one has $(u-u_0)\bot_{L^2} (u+u_0)$, so the triangle $(u_0,u,-u_0)$ is $L^2$-right at the vertex $u$, and one has
$$\sin^2\angle(\R v,\R u_0)+\sin^2\angle(\R v',\R u_0)=1\ .$$
Therefore, either $\sin\angle(\R v,\R u_0)\geq \frac{1}{\sqrt{2}}$ or $\sin\angle(\R v',\R u_0)\geq \frac{1}{\sqrt{2}}$. Suppose we are in the first case.

We get the inequality:
$$\|v\|_{L^2}^2\leq 2\|{\rm pr}_{\R u_0^\bot}v\|^2\leq \frac{2}{\lambda(A_0)}\|d_A(v)\|^2\ ,
$$
where ${\rm pr}_{\R u_0^\bot}$ stands for the $L^2$-orthogonal projection on the $L^2$-orthogonal complement of the line $\R u_0=\ker(d_{A_0})$ and $\lambda(A_0)$ is the first positive eigenvalue of $\Delta_{A_0}:A^0(F)\to A^0(F)$. In other words, $\frac{1}{\sqrt{\lambda(A_0)}}$ is the norm of the inverse $G_{A_0}:\im(d_{A_0})\to \ker(d_{A_0})^\bot$ of $d_{A_0}$ with respect to the $L^2$-norms.

Using (\ref{ineq1}) we get an estimate of the form
$$E_{A_0+\alpha}(u)-E_{A_0+\alpha}(u_0)\geq \|d_{A_0}v\|^2-c_1(A_0)\sup|\alpha| \|d_{A_0}v\|^2-c_2(A_0)\sup|\alpha|^2 \|d_{A_0}v\|^2
$$
for positive constants $c_1(A_0)$, $c_2(A_0)$. Taking 
$$\varepsilon(A_0)=\min(\frac{1}{4c_1(A_0)}, \sqrt{\frac{1}{2c_2(A_0)}})\ ,$$
we get for $\sup|\alpha|\leq\varepsilon(A_0)$
$$E_{A_0+\alpha}(u)-E_{A_0+\alpha}(u_0)\geq \frac{1}{2}\|d_{A_0}v\|^2 \geq \frac{\lambda(A_0)}{4} \|v\|^2\ ,$$
which is strictly positive unless $v=0$, i.e. $u=u_0$. The same argument applies in the case  $\sin\angle(\R v',\R u_0)\geq \frac{1}{\sqrt{2}}$ by replacing $u$ with $-u_0$ and $v$ with $v'$.
\qed

The inequality $\|\alpha\|_{L^\infty}<\varepsilon(A_0)$ as $A_0$ varies in ${\cal A}'(F)$ defines a gauge invariant neighborhood ${\cal N}$ of the zero section  in the normal bundle $N$ of this submanifold.
\begin{co}\label{inj} The restriction of the natural  map 
$$\nu:N\map {\cal A}(F)\ ,\ \nu(A_0,\alpha):=A_0+\alpha$$
 to ${\cal N}$ is injective.
\end{co}
\pf Indeed, if $\nu(A_0,\alpha)=\nu(B_0,\beta)=A$, then one must have $\ker(d_{A_0})=\ker(d_{B_0})$, because the  absolute minimum of $E_A$  on  $\Gamma(S(F))$ is unique up to sign. 
By Remark \ref{A0}, we obtain $A_0=B_0$, so finally $\alpha=\beta$.
\qed

Combining with Lemma \ref{open}, we get

\begin{thry}\label{neighb1} There exists a gauge invariant neighborhood ${\cal U}$ of the zero section in   the bundle $N$   which (after suitable Sobolev completions) is mapped diffeomorphically onto its image via the natural map $\nu$.
\end{thry}
\begin{re} The map $\resto{\nu}{\cal U}$ can be regarded  as a system  of "gauge equivariant polar coordinates" around the submanifold ${\cal A}'(F)$.
\end{re}
 \begin{re} In the finite dimensional framework, one can prove easily that in general, for any submanifold $X$ of a Riemannian manifold $Y$ there exists a neighborhood of the zero section in the normal bundle $N_{X/Y}$ which is mapped diffeomorphically  via the exponential map onto a normal neighborhood of $X$.  However, in the infinite dimensional framework, the problem is much more difficult. Theorem  \ref{neighb1} solves this problem in the special case of the embedding ${\cal A}'(F)\subset {\cal A}(F)$. 
\end{re}
\subsection{Normal neighborhoods  of  the reduction loci in the space of Hermitian connections}\label{nnrl}

Let $E$ be a rank 2 Hermitian bundle over a 4 manifold $M$  and denote $D:=\det(E)$, $d:=c_1(D)$.  Consider the involution $i_d: x\mapsto d-x$ on $H^2(M,\Z)$. A congruence class $\lambda\in H^2(M,\Z)/\langle i_d\rangle$ will be called a topological decomposition of $E$ if it coincides with the set of  Chern classes of the terms of a splitting of $E$ as direct sum of line bundles (i.e. when $x(d-x)=c_2(E)$ for $x\in\lambda$).

Fix a connection $a\in{\cal A}(D)$ and denote by ${\cal A}_a(E)$ the affine space of connections on $E$ inducing $a$ on $D$. Our gauge group is the group ${\cal G}_E:=\Gamma(SU(E))$ of determinant 1 unitary isomorphisms of $E$.  A connection will be called  simply reducible if $\dim(\ker(d_A:A^0(su(E))\to A^1(su(E))))=1$. Such a connection admits precisely two parallel line subbundles (which, of course, might be isomorphic), and these subbundles give an $A$-parallel orthogonal splitting of $E$. If $A$ admits a parallel line bundle $L$ with $2c_1(L)\ne d$, then $A$ is automatically simply reducible and $L$, $L^\bot$ are the  unique $A$-parallel line subbundles of $E$. In particular, if $d\not\in 2 H^2(M,\Z)$, then any reducible connection on $E$ is simply reducible. If we fix a  line subbundle $L\hookrightarrow E$, there exists a natural bijection between the simply reducible connections of $E$  for which $L$ is parallel and the subspace ${\cal A}^*(L)\subset {\cal A}(L)$  of abelian connections $b\in {\cal A}(L)$ for which $b\not\simeq a\otimes b^{\vee}$. On has automatically ${\cal A}^*(L)= {\cal A}(L)$ when $2c_1(L)\ne d$, whereas  ${\cal A}^*(L)$ is the complement of $2^{b_1(N)}$ gauge orbits  in ${\cal A}(L)$ when $2c_1(L)= d$.

 Let $\lambda$ be a topological decomposition of $E$ and denote by ${\cal A}_a^\lambda(E)$   the subspace of  simple reducible connections $A\in {\cal A}_a(E)$ with the property that the set of the Chern classes of the two  $A$-parallel line subbundles of $E$  coincides with $\lambda$. Such a connection will be called $\lambda$-reducible. Denote also by $\Gamma^\lambda(S(su(E)))$ the set of sections $u\in \Gamma(S(su(E)))$ with the property that the set $c(u)$ of Chern classes of the eigen line subbundles of $u$ coincides with $\lambda$. Putting
\begin{equation}\label{eig}
L_u:=\ker(u-\frac{i}{\sqrt{2}}\id_E)
 \end{equation}
 these eigen line subbundles are $L_{\pm u}$.
 One has a natural surjection
 $$w:{\cal A}_a^\lambda(E)\map \qmod{\Gamma^\lambda(S(su(E)))}{\{\pm 1\}}
 $$
which associates to every connection $A\in {\cal A}_a^\lambda(E)$ the unordered pair of sections  of the sphere bundle $S(su(E))$ which are $A$-parallel.

For sufficiently large Sobolev index $k$, the space   $\Gamma(S(su(E)))_{k+1}$ becomes a Banach manifold and  the subset $\Gamma^\lambda(S(su(E)))_{k+1}$ is open and closed in $\Gamma(S(su(E)))_{k+1}$, so it is a union of connected components. The gauge group ${\cal G}_{E,k+1}$ acts smoothly on $\Gamma(S(su(E)))_{k+1}$, leaving invariant $\Gamma^\lambda(S(su(E)))_{k+1}$; on the quotient $\Gamma^\lambda(S(su(E)))_{k+1}/\{\pm 1\}$ this gauge group  acts transitively.
The same arguments as in the proof of Proposition \ref{subm} show that
\begin{pr} \label{subm1}   The subset  ${\cal A}_a^\lambda(E)_k\subset{\cal A}_a(E)_k$ is a submanifold, which is a locally trivial affine bundle over  $\Gamma^\lambda(S(su(E)))_{k+1}/\{\pm 1\}$. The fibre over a class $[u]\in \Gamma(S(su(E)))_{k+1}/\{\pm 1\}$ can be naturally identified with the space of  abelian connections ${\cal A}^*(L_u)_k$.
\end{pr}
 Fix a reducible connection $A\in{\cal A}_a^\lambda(E)$ with $w(A)=[u]$. Put $S_u:=(L^{\otimes 2}_u\otimes D^\vee)$.
The bundle $su(E)$ splits as an orthogonal sum   sum of $A$-parallel summands:
$$su(E)=(M\times \R u)\oplus S_u=so(u^\bot)\oplus so(u^\bot)^\bot\ ,$$
and, as in the  section \ref{nn}, we obtain the following $L^2$-orthogonal decomposition of  the tangent space $T_A({\cal A}_a(E))=A^1(su(E))$ at   $A$:
$$A^1(su(E))=A^1(M,\R) u\oplus A^1(S_u)=$$
$$A^1(M,\R) u\oplus d_A [A^0(S_u)]\oplus \ker \left[d^*_A:A^1(S_u)\to A^0(S_u)\right].
$$ 
The first summand $A^1(M,\R)u$ is the tangent space of the fibre $w^{-1}(u)\simeq{\cal A}^*(L_u)$ at $A$, $n_A:=d_A[A^0(S_u)]$ is the normal space at $A$ of this fibre in the submanifold ${\cal A}_a^\lambda(E)$, whereas the space $N_A:=\ker\left[d^*_A:A^1(S_u)\to A^0(S_u)\right]$ is the normal space at $A$  of  ${\cal A}_a^\lambda(E)$ in the space of connections ${\cal A}_a(E)$.  %

We denote by $N^\lambda \to{\cal A}_a^\lambda(E)$ the $L^2$-normal vector bundle of the submanifold ${\cal A}_a^\lambda(E)$, whose fibre over $A\in {\cal A}_a^\lambda(E)$ is  $N_A$.\\
 
 Using Theorem  \ref{neighb1} we obtain the following important result, which gives an $L^2$-normal  neighborhood of the submanifold ${\cal A}_a^\lambda(E)$ of simple reductions  of type $\lambda$, and a system of polar coordinates around this submanifold.
\begin{thry}\label{neighb}   There exists a gauge invariant neighborhood ${\cal U}^\lambda$ of the zero section in   the normal bundle $N^\lambda\to{\cal A}_a^\lambda(E)$   which (after suitable Sobolev completions) is mapped diffeomorphically onto its image via the natural map 
$$\nu:N^\lambda\to {\cal A}_a(E)\ ,\ (A,\alpha)\mapsto A+\alpha\ .$$

The elements of the orthogonal slice $\nu(N_A\cap{\cal U}^\lambda)\subset{\cal A}_a(E)$ through $A\in  {\cal A}_a^\lambda(E)$ are connections for which the two elements of 
$$\ker(d_A:A^0(su(E))\to A^1(su(E)))\cap \Gamma(S(su(E)))$$
 (which are $A$-parallel) are the unique energy-minimizing harmonic sections in $\Gamma(S(su(E)))$.
\end{thry}

From now on we will always assume that ${\cal U}^\lambda$ is defined by an inequality of the form $\|\alpha\|_{L^\infty}\leq \varepsilon(A)$, where $A\mapsto \varepsilon(A)$ is continuous and gauge invariant  (see Lemma \ref{open} and Lemma \ref{absmin}).

\subsection{Twisted reductions} \label{twred}

When the base manifold has nontrivial first homology group $H_1(M,\Z)$, one also has to take into account the twisted reductions, i.e. the connections which are irreducible but whose pull-back on a double cover of $M$  become reducible. Although the stabilizer ${\cal G}_{E,A}$ of such a connection is just the center $\{\pm \id_E\}$ of the gauge group, these loci of twisted connections and the geometry of the instanton moduli spaces  around these loci must be studied in detail; the reason is  simple: the classical  transversality results with respect to metric variations \cite{DK} fails not only at a reduction, but also at a twisted reduction, so it is not clear whether one can achieve regularity of an instanton moduli space at such a point by perturbing the metric.

Let $\rho:\pi_1(M,x_0)\to \{\pm 1\}$ be a group epimorphism and denote by $\pi_\rho:M_\rho\to M$ the double cover associated with $\ker(\rho)$. The tautological involution of $M_\rho$ will be denoted by $\iota$.
A connection $A\in{\cal A}_a(E)$ will be called {\it $\rho$-twisted  reducible} (or a $\rho$-{\it twisted reduction}) if it is irreducible, but its pullback $\pi_\rho^*(A)\in {\cal A}_{\pi_\rho^*(a)}(\pi_\rho^*(E))$ is reducible.  For such a connection one has an orthogonal $A$-parallel splitting $\pi_\rho^*(E)=L'\oplus L''$ and an isomorphism $L''\simeq \iota^*(L')$.  $L'$ and $L''$ are the unique $\pi_\rho^*(A)$-parallel subbundles of $E$, because, if not, the $SO(3)$-connection associated with $\pi_\rho^*(A)$  would be trivial. In this case $A$ must be  flat and the holonomy of $A$ acts on the projective line $\P(E_{x_0})$ by an involution. Therefore,  $A$ would admit (at least) two parallel line subbundles, contradicting the assumption that $A$ was irreducible. In other words $\pi_\rho^*(A)$ must be simply reducible of type $\lambda=(l,\iota^*(l))$, where $l$ is a solution in $H^2(M_\rho,\Z)$ of the system
\begin{equation}\label{system-l}
l+\iota^*(l)=\pi_\rho^*(d)\ ,\ l\cdot \iota^*(l)=2c_2(E)\ .
\end{equation}
For a $\rho$-twisted reduction $A$ the $SO(3)$-bundle $su(E)$ has an $A$-parallel splitting 
$$su(E)=R_\rho\oplus  F\ ,
$$ 
where $R_\rho$ is the Euclidean real  line subbundle of $su(E)$ consisting  of trace-free anti-Hermitian endomorphisms $u_x\in su(E_x)$ whose eigen lines are   $L'_{\tilde x}$, $L''_{\tilde x}=L'_{\iota(\tilde x)}$, where $\tilde x$ is a lift of $x\in M$ in $M_\rho$, and $L'$, $L''$ are the   $\pi_\rho^*(A)$-parallel line subbundles of $\pi_\rho^*(E)$. $R_\rho$ is isomorphic with the non-orientable Euclidean real  line bundle associated with the representation $\rho:\pi_1(M,x_0)\to O(1)$. The second term $F$ is an $O(2)$-bundle with $\det(F)\simeq R_\rho$, whose pull-back to $M_\rho$ has two $SO(2)=U(1)$-reductions, isomorphic with $[L']^\vee \otimes   L''$ and  $[L'']^\vee \otimes   L'$.

We fix  a topological decomposition $\lambda=\{l,\iota^*(l)\}$ of $\pi_\rho^*(E)$ and we denote by ${\cal A}^\lambda_a(E)$ the subspace of $\rho$-twisted reductions $A\in {\cal A}_a(E)$ with the property that $\pi_\rho^*(A)$ is simply reducible of type $\lambda$.

Denote by $\Gamma_\iota(S(su(\pi_\rho^*(E))))$  the set of sections $u$ of the sphere bundle of $su(\pi_\rho^*(E))$ satisfying  the property $\iota^*(u)=-u$ and by $\Gamma_\iota^\lambda(S(su(\pi_\rho^*(E))))$ the subset of $\Gamma_\iota(S(su(\pi_\rho^*(E))))$ consisting of sections $u$, such that the  Chern classes of the eigen line sub-bundles of $\pi_\rho^*(u)$ are $l$, $\iota^*(l)$.

In the same way as in the case of non-twisted connections one gets a locally trivial, gauge equivariant fibration
$$
w:{\cal A}^\lambda_a(E)\map \qmod{\Gamma_\iota^\lambda(S(su(\pi_\rho^*(E))))}{\{\pm 1\}}\ .
$$
and the fibre over a class $[u]$ can be identified with the subspace 
$${\cal A}_\iota^*(L_{\pi_\rho^*(u)})\subset {\cal A}(L_{\pi_\rho^*(u)})$$
of abelian connections $b\in {\cal A}(L_{\pi_\rho^*(u)})$ satisfying 
$$b\otimes \iota^*(b)=\iota^*(a)\hbox{ (via the obvious isomorphism $L\otimes \iota^*(L)\simeq \pi_\rho^*(D)$) },\ b\not\simeq \iota^*(b).$$
The second condition is superfluous when $l\ne\iota^*(l)$. The   space $A^1(su(E))$ splits as
$$A^1(su(E))=A^1(R_\rho)\oplus d_A [A^0(F)]\oplus \ker \left[d^*_A:A^1(F)\to A^0(F)\right]\ ,
$$
where the third term can be identified with the normal space at $A$ of the submanifold ${\cal A}^\lambda_a(E)$ of $\rho$-twisted  connections of type $\lambda$.

Using similar methods as in the proofs of Theorems \ref{neighb1}, \ref{neighb} one gets easily the following existence theorem for $L^2$-normal neighborhoods of the loci of twisted reductions.
\begin{thry}\label{neighbtw}   There exists a gauge invariant neighborhood ${\cal U}^\lambda$ of the zero section in   the normal bundle $N^\lambda\to{\cal A}_a^\lambda(E)$   which (after suitable Sobolev completions) is mapped diffeomorphically onto its image via the natural map 
$$\nu:N^\lambda\to {\cal A}_a(E)\ ,\ (A,\alpha)\mapsto A+\alpha\ .$$

The elements of the orthogonal slice $\nu(N_A\cap{\cal U}^\lambda)\subset{\cal A}_a(E)$ through $A\in  {\cal A}_a^\lambda(E)$ are connections for which the two elements of the intersection
$$\ker(d_{\pi_\rho^*(A)}:A^0(su(\pi_\rho^*(E)))\to A^1(su(E)))\cap \Gamma_\rho(S(su(\pi_\rho^*(E))))$$
 (which are $\pi_\rho^*(A)$-parallel) are the unique energy-minimizing harmonic sections in $\Gamma(S(su(E)))$.
 \end{thry}
 
 As in the non-twisted case we will suppose that ${\cal U}^\lambda$ is defined by an inequality of the form $\|\alpha\|_{L^\infty}\leq \varepsilon(A)$, where  the assignment ${\cal A}_a^\lambda(E)\ni A\mapsto \varepsilon(A)$ is gauge invariant and continuous.

\section{The Donaldson $\mu$-classes around the reductions and the homotopy type of ${\cal B}_a(E)$}

\subsection{The universal bundle around reductions}\label{ub}

The structure of the universal bundle around a single reduction  is well-known. A complete description can be found in \cite{DK} p. 186-187. However, for our purposes, this classical result is not sufficient, because we will need the structure of the universal bundle around positive dimensional subspaces of reductions. \\

We recall that the universal $SO(3)$-bundle on ${\cal B}^*_a(E)\times M$ is defined as
$$\F:=\qmod{{\cal A}^*_a(E)\times su(E)}{\bar{\cal G}_E}\ ,
$$
where $\bar{\cal G}_E:={\cal G}_E/\{\pm 1\}$ acts in the natural way on both factors. Alternatively, one can let  $\bar{\cal G}_E$  act from the right and define $\F$  to be the bundle  with fibre $su(E)$ over ${\cal B}^*_a(E)$ which is  associated with  the principal $\bar{\cal G}_E$-bundle ${\cal A}^*_a(E)\to {\cal B}^*_a(E)$.   Let $\lambda=\{l,d-l\}$ be a topological decomposition of $E$.  We will assume for simplicity  that $2l\ne d$ which assures that ${\cal A}(L)={\cal A}^*(L)$ for every Hermitian line bundle of Chern class $l$; in particular the fibres of the fibration $w:{\cal A}^\lambda_a(E)\to \Gamma^\lambda(S(su(E))/\{\pm 1\}$ are affine spaces of the form ${\cal A}(L_u)$ (see section \ref{nnrl}).

We will omit the upper script   $\lambda$ in the notations $N^\lambda$, ${\cal U}^\lambda$ introduced in the previous section and we denote by $N^*$, ${\cal U}^*$ the complement of the zero section in $N$ and ${\cal U}$ (see Theorem \ref{neighb}).  Replacing ${\cal U}$ be a smaller gauge invariant neighborhood if necessary, we may assume that $\nu({\cal U}^*)\subset  {\cal A}^*_a(E)$.  We put 
$${\cal V}:=\qmod{\nu({\cal U})}{\bar{\cal G}_E}\subset {\cal B}_a(E)\ ,\ {\cal V}^*:=\qmod{\nu({\cal U}^*)}{\bar{\cal G}_E}\subset {\cal B}^*_a(E)\ .
$$
${\cal V}$ is a neighborhood of the moduli space  of $\lambda$-reductions  ${\cal B}_a^\lambda(E):={\cal A}^\lambda_a(E)/{\cal G}_E$.
Fix a section $u\in  \Gamma^\lambda(S(su(E)))$ and denote by ${\cal G}_u$ the subgroup  of ${\cal G}_E$ consisting of elements $g\in{\cal G}$ which leave $L_u$ invariant.  We have a natural isomorphism ${\cal G}_u\simeq{\cal G}$, where ${\cal G}:={\cal C}^\infty(M,S^1)$. Put 
$$\bar{\cal G}:=\qmod{{\cal G}}{\{\pm 1\}}\ ;\ \bar{\cal G}_u:=\qmod{{\cal G}_u}{\{\pm 1\}}\subset\bar{\cal G}_E\ ,$$
and  denote by ${\cal G}^{x_0}$, ${\cal G}_E^{x_0}$, ${\cal G}_u^{x_0}$ the subgroups of  ${\cal G}$ (respectively ${\cal G}_E$, ${\cal G}_u$)  of elements $f$ with $f(x_0)=1$. Note that these subgroups are mapped injectively into $\bar{\cal G}$, $\bar{\cal G}_E$, $\bar{\cal G}_u$, so we will use the same notations for the corresponding subgroups  of these groups.

The main point which will be used in our computation is that the fixing of the section $u$ defines a $\bar{\cal G}_u$-reduction of the restriction of the principal $\bar{\cal G}_E$-bundle ${\cal A}^*_a(E)\to {\cal B}^*_a(E)$ to the subspace ${\cal V}^*\subset {\cal B}^*_a(E)$.
\begin{pr}\label{reduction} Suppose $\lambda=\{l,d-l\}$ with $2l\ne d$. Let $u\in \Gamma^\lambda(S(su(E)))$ and let   $N_u\subset {\cal A}(L_u)\times A^1(S_u)$ the restriction of the normal bundle $N$ to the fibre $w^{-1}([u])\simeq{\cal A}(L_u)$. Then
\begin{enumerate}\label{red}
\item   $\nu(N_u\cap {\cal U})$ is the submanifold of the normal neighborhood $\nu({\cal U})$   consisting of connections for which   $\pm u$ are harmonic and energy minimizing.
\item \label{resu} The embedding $\resto{\nu}{N_u\cap {\cal U}}:N_u\cap {\cal U}\hookrightarrow {\cal A}_a(E)$ induces isomorphisms
$${\cal V}\simeq \qmod{ N_u\cap{\cal U}}{\bar{\cal G}_u}\ , {\cal V}^*\simeq \qmod{N^*_u\cap{\cal U}}{\bar{\cal G}_u}\ .
$$
\item The map $\nu$ induces an isomorphism  between ${\cal V}$ and the cone bundle over the projectivization  $\P(\bar N_u)$ of  the vector bundle 
$$\bar N_u:=\qmod{N_u}{{\cal G}_u^{x_0}}
$$
over ${\cal B}(L_u)=\qmod{{\cal A}(L_u)}{{\cal G}_u^{x_0}}$, and ${\cal V}^*$ is identified with the complement of the vertex section in this cone bundle. In particular one has a homotopy equivalence
\begin{equation}\label{he}
{\cal V}^*\stackrel{\rm h}{\simeq}\P(\bar N_u)\subset {\cal A}(L_u)\times_{{\cal G}_u^{x_0}}\P(A^1(S_u))\ .
\end{equation}
\item $\nu(N_u\cap {\cal U}^*)$ is a $\bar{\cal G}_u$-reduction of the restriction of the principal $\bar{\cal G}_E$-bundle ${\cal A}^*_a(E)\to {\cal B}^*_a(E)$ to  ${\cal V}^*$.
\end{enumerate} 
\end{pr}
\pf The first statement follows easily from Theorem \ref{neighb}. The second and the third statements are obvious. For the fourth, it suffices to prove that the $\bar{\cal G}_E$-orbit of a point $\nu(t)\in\nu({\cal U}^*)$ intersects $\nu(N_u\cap {\cal U}^*)$ along a $\bar{\cal G}_u$-orbit or, equivalently, that   the $\bar{\cal G}_E$-orbit of $t=(b,\alpha)\in N^*$ intersects $N_u^*$ along along a $\bar{\cal G}_u$-orbit. Since $\bar{\cal G}_E$ acts transitively on $\Gamma^\lambda(S(su(E)))/\{\pm 1\}$, we see that $t\cdot \bar{\cal G}_E\cap N_u^*\ne\emptyset$, so we can suppose $t\in N_u^*$.  If an element $g\in{\cal G}_E$ maps $t$ into $N_u^*$, then $g(b)\in{\cal A}(L_u)\subset {\cal A}_a(E)$, i.e.  $u$ is both $b$-harmonic and   $g(b)$-harmonic. This implies $\ad_g(u)=\pm u$. Since we assumed that $l\ne d-l$, we have $L_u\not\simeq L_{-u}$, so necessarily $\ad_g(u)=u$ i.e. $g\in{\cal G}_u$.
\qed
\begin{co}\label{ht} The open subspace ${\cal V}^*\subset {\cal B}^*_a(E)$ has the homotopy type of  the product $[S^1]^{b_1(M)}\times \P^\infty$.
\end{co} 
\pf  Let $g$ be a Riemannian metric on $M$.  The moduli space ${\cal M}^{YM}_g(L_u)$  of Yang-Mills connections on $L_u$ with respect to $g$ is isomorphic to 
$$\qmod{i H^1(X,\R)}{2\pi i H^1(M,\Z)}\simeq [S^1]^{b_1(M)}
$$
and the inclusion  ${\cal M}^{YM}_g(L_u)\hookrightarrow{\cal B}(L_u)$ is a homotopy equivalence. On the other hand, since ${\cal M}^{YM}_g(L_u)$ is compact, the restriction of the infinite rank vector bundle  $\bar N_u$ to this subspace is trivial. This completes the proof.
\qed

The decomposition 
$$su(E)=[M\times(i\R)]\oplus S_u
$$
is ${\cal G}_u$-invariant.  Therefore
\begin{co} The restriction $\resto{\F}{{\cal V}^*}$ decomposes as a direct sum 
$$\resto{\F}{{\cal V}^*}\simeq [({\cal V}^*\times M)\times i\R]\oplus {\cal S}_u\ ,
$$
where ${\cal S}_u$ is the complex line bundle on ${\cal V}^*\times M$ defined by
$${\cal S}_u:=\qmod{\nu(N_u\cap {\cal U}^*)\times S_u}{\bar{\cal G}_u}\ .
$$
In particular $p_1(\F)=c_1({\cal S}_u)^2$.
\end{co}

The cohomology algebra of ${\cal V}^*$ can be easily described explicitly using Corollary \ref{ht}. The construction below yields generators with explicit geometric interpretation.

Put $L=L_u$, $S=S_u=L^{\otimes 2}\otimes D^\vee$, ${\cal S}={\cal S}_u$ to save on notations.  The subbundle $N_u\subset{\cal A}(L)\times A^1(S)$  has an intrinsic interpretation in terms of $L$: it is just the bundle ${\cal K}=\ker \delta$ of kernels of the family of operators 
$$\delta:=(d_{b^{\otimes 2}\otimes a^\vee}^*)_{b\in{\cal A}(L)}$$
and ${\cal U}\cap N_u$ can be identified with a $\bar{\cal G}$-invariant neighborhood ${\cal U}_u$ of the zero section in this bundle. The bundle ${\cal K}$ descends to a bundle $\bar{\cal K}:={{\cal K}}/{{\cal G}^{x_0}}$ over ${\cal B}(L)$.
Using the isomorphism ${\cal V}^*\simeq \qmod{{\cal U}^*_u}{\bar {\cal G}}$ induced by $\nu$, one can identify  ${\cal S}$ with the  line bundle 
$$\qmod{{\cal U}^*_u\times S}{\bar{\cal G}}\map \qmod{{\cal U}^*_u}{\bar{\cal G}}\times M\subset \qmod{{\cal K}^*}{\bar{\cal G}}\times M=\qmod{\bar{\cal K}^*}{S^1}\times M\ .
$$
Therefore, ${\cal S}$ can be regarded as the restriction to $[{\cal U}^*_u/\bar{\cal G}]\times M$ of the line bundle
$$\S:=\qmod{{\cal A}(L)\times (A^1(S)\setminus\{0\})\times S}{\bar {\cal G}}
$$
on  ${\cal B}^*_{u}(E)\times M$, where ${\cal B}^*_{u}(E)$ is the quotient
$${\cal B}^*_{u}(E):=\qmod{{\cal A}(L)\times (A^1(S)\setminus\{0\})}{\bar {\cal G}}={\cal A}(L)\times_{{\cal G}^{x_0}}\P(A^1(S))\ ,$$
which can be regarded as a projective bundle over ${\cal B}(L)$. The space $ {\cal B}^*_u(E)$ is very much similar to the infinite dimensional gauge quotient of the space of irreducible configurations in Seiberg-Witten theory. More precisely, let in general $V$ be a line bundle and $W$ a complex vector bundle on $M$. 
 The natural map
$$p_{V,W}:{\cal B}^*(V,W):=\qmod{{\cal A}(V)\times [A^0(W)\setminus\{0\}] } {{\cal G}}\map {\cal B}(V)
$$
is a projective bundle. On the product ${\cal B}(V)\times M$ one has a tautological line bundle defined  by
$$\V_W:=\qmod{{\cal A}(V)\times [A^0(W)\setminus\{0\}]\times V } {{\cal G}}
$$

\begin{lm} With the notations and assumptions above one has
\begin{enumerate}
\item There exists a natural isomorphism 
$$\delta_V: H_1(M,\Z)\to H^1({\cal B}(V))  $$
 which induces  an  isomorphism $H^*({\cal B}^*(V),\Z)\simeq\Lambda^*(H_1(M,\Z))$.
\item The morphism $p_{V,W}^*\circ\delta_V$ induces  an  isomorphism 
$$H^*({\cal B}^*(V,W),\Z)\simeq\Lambda^*(H_1(M,\Z))\otimes\Z[h_{VW}]\ ,
$$
where $h_{VW}$ is a degree 2-cohomology class defined as the Chern class of the principal $S^1$-bundle
$$\qmod{{\cal A}(V)\times [A^0(W)\setminus\{0\}] } {{\cal G}^{x_0}}\map \qmod{{\cal A}(V)\times [A^0(W)\setminus\{0\}] } {{\cal G}}={\cal B}^*(V,W)\ .
$$
This class restricts to the canonical (tautological) class of the projective fibres of  ${\cal B}^*(V,W)$.
\item The Chern class of the line bundle $\V_W$ on ${\cal B}^*(V,W)\times M$ is
$$c_1(\V_W)=p_1^*(h_{VW})\otimes 1+\delta_V+1\otimes p_2^*(c_1(V))\ ,
$$
where we agree to denote by the same symbol $\delta_V$ the   element  in $H^1({\cal B}(V))\otimes H^1(M,\Z)$ defined by this morphism, as well as  the pullback of this element via the projection ${\cal B}^*(V,W)\times M\to {\cal B}(V)\times M$.
\end{enumerate}
\end{lm}
\pf The proof uses the same arguments which are used for the computation of the cohomology algebra of the moduli space of irreducible configurations in Seiberg-Witten theory (see for instance \cite{OT1}).
\qed
\begin{co} The line bundle $\S$ can be identified with the pull-back of $\S_{\Lambda^1\otimes S}$ via the composition
$${\cal B}^*_u(E)\map  {\cal B}^*(S,\Lambda^1\otimes S)
$$
induced by  
$${\cal A}(L)\times [A^1(S)\setminus\{0\}]\ni (b,\alpha)\mapsto (b^{\otimes  2}\otimes a^{-1},\alpha)\in {\cal A}(S)\times [A^1(S)\setminus\{0\}]\ .$$
 In particular one has
$$c_1(\S)=p_1^*(h_u)\otimes 1+ 2\delta_L+p_2^*(c_1(S))\ ,
$$
where $h_u$ is the Chern class of the principal $\bar {\cal G}/{\cal G}^{x_0}=S^1$-bundle
$$\qmod{{\cal A}(L)\times (A^1(S)\setminus\{0\})}{ {\cal G}^{x_0}}\map\qmod{{\cal A}(L)\times (A^1(S)\setminus\{0\})}{\bar {\cal G}}={\cal B}^*_u(E)\ .
$$
\end{co}
\pf  The map $(b,\alpha)\mapsto (b^{\otimes  2}\otimes a^{-1},\alpha)$ induces an  isomorphism $u$ which fits in the diagram
$${\cal B}^*_u(E)= \qmod{{\cal A}(L)\times [A^0(\Lambda^1\otimes S)\setminus\{0\}] } {\bar{\cal G}}\textmap{u} \qmod{{\cal A}(S)\times [A^0(\Lambda^1\otimes S)\setminus\{0\}] } {{\cal G}^2}\textmap{v}$$
$$\textmap{v} \qmod{{\cal A}(S)\times [A^0(\Lambda^1\otimes S)\setminus\{0\}] } {{\cal G}}={\cal B}^*(S,\Lambda^1\otimes S)\ ,$$
where $v$ is the obvious epimorphism. The map 
$${\cal A}(L)\times [A^1(S)\setminus\{0\}]\times S\in (b,\alpha, s)\mapsto (b^{\otimes  2}\otimes a^{-1},\alpha,s)\in {\cal A}(S)\times [A^1(S)\setminus\{0\}]\times S$$
 descends to an isomorphism  $\S\simeq (v\circ u)^*(\S_{\Lambda^1\otimes S})$. It suffices to  notice that $(v\circ u)^*(h_{VW})=h_u$ and  $(v\circ u)^*(\delta_S)=2\delta_L$.
\qed
\begin{co}\label{dclasses} Suppose that $b_+(M)=0$. The restrictions of the Donaldson $\mu$-classes to ${\cal V}^*$ are given  by the formulae
$$\mu(\lambda)=-  \langle \delta_L\cup_M c_1(S),\lambda\rangle\in H^1({\cal V}^*,\Z)\   \forall \lambda\in H_3(M,\Z)\ ,
$$
$$\mu(x)=-\frac{1}{2}\langle c_1(S),x\rangle  h_u \in H^2({\cal V}^*,\Q) \ \forall x\in H_2(X,\Z)\ ,
$$
$$\mu(\sigma)=-\delta_L(\sigma)\cup h_u \in H^3({\cal V}^*,\Z)\ \forall \sigma\in H_1(X,\Z)\ ,
$$
$$\nu:=\mu(*)=-\frac{1}{4}   h_u^2\in H^4({\cal V}^*,\Q)\ .
$$
\end{co}
\pf Every monomial of the form $a\cup b$, with $a$, $b\in H^1(B,\Z)$ satisfies the relation $(a\cup b)^2=0$  so,  since we supposed $b_+(M)=0$, one gets easily that $a\cup b\in\Tors(H^2(M,\Z))$. Therefore 
$$\im(\cup:H^1(M,\Z)\otimes H^1(M,\Z)\to  H^2(M,\Z))\subset\Tors(H^2(M,\Z)).$$
  This implies $\delta_L^2=0$. Therefore
$$-\frac{1}{4}p_1(\resto{\F}{{\cal V}^*})=-\frac{1}{4}c_1(\S)^2=$$
$$-\frac{1}{4}[p_1^*(h_u^2)\otimes 1+4p_1^*(h_u)\cup \delta_L+2p_1^*(h_u)\cup p_2^*(c_1(S))+ 4\delta_L\cup p_2^*(c_1(S))]\ .
$$
Using Donaldson's formula $\mu(x)=-\frac{1}{4}p_1(\F)/x$ for $x\in H_*(X,\Z)$, one gets easily the claimed formulae.
\qed

\subsection{The topology of the moduli space ${\cal B}_a(E)$} \label{topology}

Describing the weak homotopy type of the moduli space ${\cal B}_a^*(E)$ of {\it irreducible} connections  is a well-known  classical problem in gauge theory. This problem is treated in detail in \cite{DK}, where the authors also compute the rational cohomology algebra of this space. Surprisingly, describing the weak homotopy type  of the whole moduli space ${\cal B}_a(E)$ of connections is a  delicate problem, which, to our knowledge, cannot be solved with similar methods.

Our result concerning the existence of $L^2$-normal neighborhoods of the reduction loci, gives a solution to this problem. Suppose for  simplicity that $d:=c_1(D)$ is not divisible by 2 in $H^2(M,\Z)$, so that any reducible connection $A\in {\cal A}_a(L)$ is simply reducible.
Let $\Lambda_E$ be the set topological decomposition of $E$, i.e. the set of unordered pairs $\lambda=\{l,d-l\}$ with $l(d-l)=c_2(E)$.    

Using the notations and conventions of section  \ref{ub}, we get the following decomposition of the space ${\cal B}_a(E)$ as a fibred sum:
\begin{equation}\label{diagram}
\begin{array}{ccccc}
&&{\cal B}_a(E)&&
\\ \\
&\nearrow&&\nwarrow&\\  \\
{\cal B}^*_a(E) && &&\hspace{-5mm}\left[\coprod_{\lambda\in \Lambda_E}{\cal V}^\lambda\right]   \ , \\  \\
&\nwarrow&&\nearrow&\\  \\
&&\left[\coprod_{\lambda\in \Lambda_E}{{\cal V}^\lambda}^*\right] &&
\end{array}
\end{equation}
where ${\cal V}^\lambda:=\qmod{{\cal U}^\lambda}{{\cal G}_E}$ is the gauge quotient of the gauge invariant $L^2$-normal  neighborhood ${\cal U}^\lambda$ of the reduction locus ${\cal A}_a^\lambda(E)$ constructed in section \ref{nnrl}.  Since the homotopy types of the terms are known (see Corollary \ref{ht}). Choosing a Hermitian  line bundle $L$ with $c_1(L)\in\lambda$, one has
$${\cal V}^\lambda\stackrel{\rm h}{\simeq}{\cal B}(L)\stackrel{\rm h}{\simeq} [S^1]^{b_1(M)}\ \ ,\ \ {{\cal V}^\lambda}^*\stackrel{\rm h}{\simeq} [S^1]^{b_1(M)}\times \P^\infty\ ,
$$
this description determines  the homotopy type of the space ${\cal B}_a(E)$. In particular one can compute the cohomology of this space using the Mayer-Vietors exact sequence.
\\ \\
{\bf Example 1.} Suppose that $H_1(M,\Z)\simeq\Z$, $b_2(M)=b_2^-(M)=2$. Under these assumptions, by Donaldson's first theorem, one can find a basis $(e_1,e_2)$ in  $H^2(M,\Z)$ such that $e_i^2=-1$, $e_1\cup e_2=0$. It will follow that $d:=e_1+e_2$ is an integral lift of the Stiefel-Whitney class $w_2(M)$. Let now $E$ be a rank 2-Hermitian bundle with $c_1(E)=d$, $c_2(E)=0$ and put again $D=\det(E)$.

Our problem is to compute the degree $k$-cohomology  of the space ${\cal B}_a(E)$, for $1\leq k\leq 4$ using the Mayer-Vietoris sequence applied to the decomposition \ref{diagram}.
Denote by $L_c$ a Hermitian line bundle with Chern class $c$. The set $\Lambda_E$ has two elements:
$$\lambda_0:=\{0,d\}\ ,\  \lambda_1:=\{e_1,e_2\} \ .
$$
Denote by ${\cal V}_0$, ${\cal V}_1$, ${\cal V}_0^*$, ${\cal V}_1^*$ the corresponding subspaces of ${\cal B}_a(E)$, and put ${\cal V}:={\cal V}_0\cup{\cal V}_1$,   ${\cal V}^*:={\cal V}_0^*\cup{\cal V}_1^*$. We get exact sequences:
$$\dots\map H^{i-1}({\cal B}_a^*(E))\oplus H^{i-1}({\cal V})\map H^{i-1}({\cal V}^*) \map H^i({\cal B}_a(E))\map $$
$$\map H^i({\cal B}_a^*(E))\oplus H^i({\cal V})\map H^i({\cal V}^*) \map\dots
$$
Using the standard description of the cohomology of ${\cal B}^*_a(E)$ (see \cite{DK}), we obtain
$$H^1({\cal B}^*_a(E))=\langle\mu(\lambda_0)\rangle_\Q\simeq\Q\ ,\ H^2({\cal B}^*_a(E))=\langle \mu(f_1),\mu(f_2)\rangle_\Q\simeq \Q^2\ ,$$
where $\lambda_0$ is a generator of $H_3(M,\Z)$ and $f_i:={\rm PD}(e_i)$. Since we assumed $b_1(X)=1$, we obtain easily using Poincar\'e duality that   
$$\im(\cup:H^1(M)\otimes  H^2(M)\to H^3(M))=0$$
 hence, by Corollary \ref{dclasses}, we obtain $\resto{\mu(\lambda_0)}{{\cal V}^*}=0$.   Note also that the restriction morphism $H^i({\cal V})\to H^i({\cal V}^*)$ is an isomorphism  for $i=0$, $1$ and is injective for all $i$. The exact sequence above for $i=1$ yields 
$$H^1({\cal B}_a(E))=\ker (H^1({\cal B}_a^*(E))\oplus H^i({\cal V})\to H^1({\cal V}^*))=\langle (\mu(\lambda_0),0)\rangle_\Q\simeq \Q\ .
$$
Therefore, the cohomology class defined by the Chern-Simons functional \cite{DK} associated with a hypersurface representing $\lambda_0$ extends to  the whole moduli space ${\cal B}_a(E)$. This  is a general phenomenon. For $i=2$, we obtain
  $$H^2({\cal B}_a(E))=\ker(H^2({\cal B}^*_a(E))\to H^2({\cal V}^*))\ .
  $$
  On the other hand by Corollary \ref{dclasses}
  $$\resto{\mu(f_i)}{{\cal V}^*_0}=\frac{1}{2}\langle e_1+e_2,f_i\rangle h_{u_0} \ ,\ \resto{\mu(f_i)}{{\cal V}^*_1}=\frac{1}{2}\langle e_2-e_1,f_i\rangle h_{u_1}\ ,
  $$
where $u_i\in \Gamma(S(su(E)))$ satisfies $L_{u_i}\simeq L_i$. Therefore $H^2({\cal B}_a(E),\Q)=0$. Taking into account that the restriction morphism $H^{2}({\cal B}_a^*(E))\to H^{2}({\cal V}^*)$ is surjective and $H^3({\cal V})=0$, the same exact sequence for $i=3$ yields:
$$H^3({\cal B}_a(E))=\ker\left[H^3({\cal B}_a^*(E))\to H^3({\cal V}^*)=H^3({\cal V}^*_0)\oplus H^3({\cal V}^*_1)\right]\ .
$$
The space $H^3({\cal V}^*_i)$ is generated by $\delta_{L_i}(\sigma_0)\cup h_{u_i}$, where $\sigma_0$ is a generator of $H_1(M,\Z)$, whereas 
$$H^3({\cal B}_a^*(E))=\langle \mu(\sigma_0),\mu(\lambda_0)\mu(f_1),\mu(\lambda_0)\mu(f_2)\rangle_\Q\ .$$
 By Corollary \ref{dclasses} one has $\resto{\mu(\sigma_0)}{{\cal V}^*_i}=-\delta_{L_i}(\sigma_0) h_{u_i}$ and we have seen that $\resto{\mu(\lambda_0)}{{\cal V}^*}=0$. 
This shows that 
$$H^3({\cal B}_a(E))=\langle  \mu(\lambda_0)\mu(f_1),\mu(\lambda_0)\mu(f_2)\rangle_\Q\simeq\Q^2\ .$$
The cohomology space $H^4({\cal B}_a(E))$  fits in the exact sequence
$$0\to \frac{H^3({\cal V}^*_0)\oplus  H^3({\cal V}^*_1)}{\langle(\delta_{L_0}(\sigma_0) h_{u_0},\delta_{L_1}(\sigma_0) h_{u_1} ) \rangle}\to H^4({\cal B}_a(E))\to  H^4({\cal B}_a^*(E)) \to H^4({\cal V}^*)
$$
The space $H^4({\cal B}_a^*(E))$ is freely generated by the five  classes  
$$\mu(*)\ ,\ \mu(f_1)^2\ ,\ \mu(f_2)^2\ ,\ \mu(f_1)\mu(f_2)\ ,  \mu(\lambda_0)\mu(\sigma_0)\ ,$$
 whereas $H^4({\cal V}^*)$ is freely generated by $h_{u_0}^2$ and $h_{u_1}^2$. This shows that
$$\ker(H^4({\cal B}_a^*(E)) \to H^4({\cal V}^*))=\langle \mu(*)+\mu(f_1)^2, \mu(*)+            \mu(f_2)^2,\ \mu(\lambda_0)\mu(\sigma_0)\rangle_\Q\simeq \Q^3\ .
$$
The quotient on the left is 1-dimensional, so $H^4({\cal B}_a(E))\simeq \Q^4$.
 \begin{re} A similar method can be used to compute the cohomology of the pair $({\cal B}^*_a(E), {\cal V}^*)$, where ${\cal V}=\coprod_{[L]\in R_E}{\cal V}^L$ is a normal neighborhood of the reducible locus in ${\cal B}_a(E)$.
 \end{re}
 {\bf Example 2.} Let $M$ be a 4-manifold with the topological properties considered in the example above. The exact sequence of the pair $({\cal B}^*_a(E), {\cal V}^*)$
$$\to H^{i-1}({\cal B}_a^*(E))\to H^{i-1}({\cal V}^*)  \to H^i({\cal B}_a^*(E),{\cal V}^*)\to H^i({\cal B}_a^*(E))\map H^i({\cal V}^*) \to 
$$
 written for  $i=4$ shows that  the natural morphism $H^4({\cal B}^*_a(E), {\cal V}^*)\to H^4({\cal B}_a(E))$ (induced by the restriction morphism $H^*({\cal B}_a(E), {\cal V})\to H^*({\cal B}_a(E))$ and the excision isomorphism $H^*({\cal B}_a(E), {\cal V})\to H^*({\cal B}^*_a(E), {\cal V}^*$) is an isomorphism.

\section{The instanton moduli space  around the reductions}

We denote by ${\cal M}_a^\ASD(E)\subset{\cal B}_a(E)$ the moduli of projectively ASD $a$-oriented connections in $E$, i.e. the moduli space 
$${\cal M}^\ASD_a(E):=\qmod{{\cal A}^\ASD_a(E)}{{\cal G}_E}\ ,\ {\cal A}^\ASD_a(E):=\{A\in{\cal A}_a(E)|\ (F_A^0)^+=0\}\ .
$$
In the first subsection we will study the intersection of this moduli space with a normal neighborhood ${\cal V}^\lambda$ of the reduction locus ${\cal B}^\lambda_a(E):={\cal A}_a^\lambda(E)/{\cal G}_E$. We will see that, in a neighborhood of a reduction locus, the instanton moduli problem is equivalent to an {\it abelian} moduli problem, which is very much similar to the Seiberg-Witten moduli problem.  We will denote by ${\cal M}_a^\lambda(E)$ the subspace ${\cal M}_a^\lambda(E):= {\cal M}_a^\ASD(E)\cap {\cal B}^\lambda_a(E)$ of $\lambda$-reducible instantons.

\subsection{An abelian gauge theoretical problem}\label{amp}

Let $L$ be a Hermitian line bundle which is isomorphic to a line subbundle of $E$, and put $S:=L^{\otimes 2}\otimes D^\vee$. Consider the moduli space ${\cal M}_a(L)$ of solutions of the  system
\begin{equation}\label{abelian}
\left\{
\begin{array}{ccc}
F_b^+-\frac{1}{2}F_a^+&=&(\alpha\wedge\bar\alpha)^+\ .\vspace{2mm}  \\  
(d_{b^{\otimes 2}\otimes a^\vee}^*,d_{b^{\otimes 2}\otimes a^\vee}^+)\alpha&=&0
\end{array}
\right.
\end{equation}
for pairs $(b,\alpha)\in{\cal A}(L)\times A^1(S)$, modulo the abelian gauge group ${\cal G}={\cal C}^\infty(M,S^1)$. This system is very much similar to the Seiberg-Witten system; indeed, the left hand operator in the first equation is elliptic and can be written as coupled Dirac operator.  The main difference is that in general there is no a priori  bound for the $\alpha$-component on the space of solutions of this system and, in general, the moduli space is not compact.

We denote by  ${\cal M}_a^{\rm red}(L)$,  ${\cal M}_a^*(L)$ the subspaces of reducible (respectively irreducible)  solutions. As in Seiberg-Witten theory ``reducible pair" means ``pair  with trivial $\alpha$-component". Therefore one has a natural identification 
 $${\cal M}^{\rm red}_a(L)\simeq T_a(L):=\qmod{{\cal A}_a(L)}{{\cal G}}\ ,  \hbox{ where }\ {\cal A}_a(L):=\left\{b\in{\cal A}(L)|\ F_b^+=\frac{1}{2} F_a^+\right\}\ ,$$
The space $T_a(L)$   is either empty, or a  $b_1(M)$-dimensional torus (when the harmonic representative of $c_1(S)=2c_1(L)-c_1(D)$ is ASD).

Put $l:=c_1(L)$ and $\lambda=\{l,d-l\}$ and fix an isomorphism $E=L\oplus (D\otimes L^\vee)$. The map 
$$\psi: (b,\alpha)\mapsto A_{b,\alpha}:=\left(\matrix{d_b&\alpha\cr -\bar\alpha&d_{a\otimes b^\vee}}\right)\in{\cal A}_a(E)$$
 descends to a map ${\cal M}_a(L)\to {\cal M}_a^\ASD(E)$.   The image of this map is the subspace consisting of those instantons which can be brought in Coulomb gauge with respect to a  $\lambda$-reducible connection. 
 Suppose again that $2l\ne d$ (such that any  reduction having a parallel line subbundle of Chern class $l$ is simple), and consider   the continuous, gauge invariant function $\varepsilon:{\cal A}^\lambda_a(E)\to \R_{>0}$ defining the normal neighborhood ${\cal U}^\lambda$ of ${\cal A}^\lambda_a(E)$  (see section \ref{nnrl}). For $b\in {\cal A}(L)$ put
 $$\varepsilon(b):=\varepsilon(A_b)\ , \hbox{ where } A_b:=b\oplus (a\otimes b^\vee)\in {\cal A}^\lambda_a(E)\ ,
 $$
 and denote ${\cal W}^\lambda$ the ${\cal G}$-invariant subspace of the configuration space ${\cal A}(L)\times A^1(S)$ defined by the inequality $\|\alpha\|_{L^\infty}\leq \varepsilon(b)$, and by ${\cal Z}^\lambda$ its ${\cal G}$-quotient.

\begin{pr}\label{equiv} The restriction of the instanton moduli problem to the normal neighborhood $\nu({\cal U}^\lambda)$ of ${\cal A}^\lambda_a(E)$  and the restriction of the abelian moduli problem (\ref{abelian}) to the neighborhood ${\cal W}^\lambda$ of ${\cal A}(L)\times\{0\}$ are equivalent moduli problems.

In particular, the map induced by $\psi$   applies ${\cal M}_a(L)\cap {\cal Z}^\lambda$ isomorphically onto the neighborhood ${\cal M}_a^\ASD(E)\cap {\cal V}^\lambda$ of  the subspace  ${\cal M}_a^\lambda(E)$ of $\lambda$-reducible instantons, and induces an isomorphism ${\cal M}^{\rm red}_a(L)\simeq {\cal M}_a^\lambda(E)$.
\end{pr}  
 \pf This follows directly from Proposition \ref{red}, \ref{resu}.
 \qed
\begin{co}  Let $b\in {\cal A}_a(L)$. The deformation elliptic complex ${\cal C}_{A_b}$ at the corresponding  reduction $A_b$ splits as a direct sum  ${\cal C}_{A_b}={\cal C}_0^+(M)\oplus {\cal C}_{b^{\otimes 2}\otimes a^\vee}^+(S)$, where ${\cal C}_0^+(M)$ is the standard $d^+$-elliptic complex for $i\R$-valued forms and ${\cal C}_{b^{\otimes 2}\otimes a^\vee}^+(S)$ is the $d^+$-elliptic complex for $S$-valued forms associated with the connection ${b^{\otimes 2}\otimes a^\vee}$.
\end{co}
\begin{re} The complex index of the elliptic complex ${\cal C}_{b^{\otimes 2}\otimes a^\vee}^+(S)$ is given by
$${\rm ind}_\C(C^+(S))=c_1^2(S)+(b_+(M)-b_1(M)+1)\ .
$$
\end{re}
\subsection{A twisted abelian gauge theoretical problem}\label{tamp}

Let $\rho:\pi_1(M,x_0)\to \Z_2$ be a group epimorphism and  consider the associated objects $\pi_\rho:M_\rho\to M$, $\iota:M_\rho\to M_\rho$ introduced in   section \ref{twred}.

Let $L$ be a Hermitian line bundle on $M_\rho$ whose Chern class $l$ is a solution  of the system (\ref{system-l}) and put $\lambda=\{l,\iota^*(l)\}$.  In this section we will assume that  $l\ne \iota^*(l)$. 
We fix an isomorphism $\pi_\rho^*(E)=L\oplus \iota^*(L)$. We  denote by ${\cal A}_\iota(L)$ the subspace of  ${\cal A}(L)$ consisting of connections $b\in {\cal A}(L)$  such that $b\otimes \iota^*(b)=\pi_\rho^*(a)$. The natural gauge group acting on this space of connections is
$${\cal G}_\iota:=\{f\in {\cal C}^\infty(M_\rho,S^1)|\ \iota^*(f)= \bar f\}\ .
$$
The Lie algebra of this  group can be identified with $A^0(R_\rho)$, where $R_\rho$ is the real line bundle associated with the representation $\rho$. One has a natural embedding ${\cal G}_\iota\to {\cal G}_{\pi_\rho^*(E)}$ factorizes   through an embedding   ${\cal G}_\iota\hookrightarrow {\cal G}_E$.  

One has a  ${\cal G}_\iota$-equivariant map ${\cal A}_\iota(L)\to {\cal A}^\lambda_a(E)$ given by $b\mapsto A_b$, where $A_b$ is the unique $\rho$-twisted connection of type $\lambda$ whose pull-back to $M_\rho$ is $b\oplus \iota^*(b)$.

Put $S=L\otimes\iota^*(L)^\vee=L^{\otimes 2}\otimes\pi_\rho^*(D)$. This Hermitian line bundle comes with a tautological isomorphism $\iota^*(S)=S^\vee=\bar S$. We introduce the spaces of $\iota$-twisted $S$-valued forms  by
$$A^k_\iota(S):=\{\alpha\in A^k(S)|\ \iota^*(\alpha)=-\bar \alpha\}\ .$$
Note that, for a connection $b\in {\cal A}_\iota(L)$, one has $F_b-\frac{1}{2}\pi_\rho^*(F_a)\in A^2_\iota(S)$.
Endow $M_\rho$ with the pull-back metric $\pi_\rho^*(g)$. Our abelian moduli problem is now
\begin{equation}\label{twabelian}
\left\{
\begin{array}{ccc}
F_b^+-\frac{1}{2}\pi_\rho^*(F_a)^+&=&(\alpha\wedge\bar\alpha)^+\ 
\vspace{2mm}  \\  
(d_{b^{\otimes 2}\otimes \pi_\rho^*(a)^\vee}^*,d_{b^{\otimes 2}\otimes \pi_\rho^*(a)^\vee}^+)\alpha&=&0\ .
\end{array}
\right.
\end{equation}
for pairs $(b,\alpha)\in{\cal A}_\iota(L)\times A^1_\iota(S)$, modulo the abelian gauge group ${\cal G}_\iota$.

The stabilizer of {\it any} pair $(b,\alpha)\in{\cal A}_\iota(L)\times A^1_\iota(S)$ is $\{\pm 1\}\subset {\cal G}_\iota$, so, in fact, this moduli twisted abelian problem has no reductions at all. However, the moduli space ${\cal M}_{a,\iota}(L)$ of solutions has two distinguished gauge invariant  subspaces ${\cal M}_{a,\iota}^{\rm red}(L)$, ${\cal M}_{a,\iota}^*(L)$ consisting of classes of solutions with vanishing (respectively non-vanishing) $\alpha$-component. As in the previous section we obtain a map $\psi: {\cal M}_{a,\iota}(L)\to {\cal M}^{ASD}_a(E)$ which applies isomorphically  ${\cal M}_{a,\iota}^{\rm red}(L)$ onto the space  ${\cal M}^{\lambda}_a(E)$ of $\rho$-twisted reducible instantons  of type $\lambda$. 

The space ${\cal M}_{a,\iota}^{\rm red}(L)$ has a simple geometric interpretation: one has an obvious identification
$${\cal M}_{a,\iota}^{\rm red}(L)\simeq \qmod{{\cal A}_{\pi_\rho^*(a)}(L)\cap {\cal A}_\iota(L)}{{\cal G}_\iota}\subset {\cal M}_{\pi_\rho^*(a)}^{\rm red}(L)
$$
which shows that  ${\cal M}_{a,\iota}^{\rm red}(L)$ is the subspace of   ${\cal M}_{\pi_\rho^*(a)}^{\rm red}(L)$ defined by the equation $[b\otimes\iota^*(b)]=[a]$. If  the harmonic representative of  $c_1(S)$ is not ASD (with respect to the pull-back metric), the space  ${\cal M}_{a,\iota}^{\rm red}(L)$ will be empty.  When this  representative is ASD,   ${\cal M}_{\pi_\rho^*(a)}^{\rm red}(L)$ will be a (non-empty!) subtorus of the torus ${\cal M}_{\pi_\rho^*(a)}^{\rm red}(L)$; this subtorus  is (non-canonically) isomorphic to the quotient
$$\qmod{ iH^1(M_\rho,\R)_\iota}{2\pi iH^1(M_\rho,\Z)_\iota}\ ,$$
 where  the symbol  $(-)_\iota$ means $\iota$-twisted, i.e. the  subspace of $(-)$ consisting of solutions of the equation $\iota^*(x)=-x$. In particular, if $\iota^*:H^1(M_\rho,\Z)\to H^1(M_\rho,\Z)$ is the identity and $b_+(M_\rho)=0$, this space consists of a single point.

Consider  the continuous, gauge invariant function $\varepsilon:{\cal A}^\lambda_a(E)\to \R_{>0}$ defining the normal neighborhood ${\cal U}^\lambda$ (see section \ref{twred}) of ${\cal A}^\lambda_a(E)$ and,  for $b\in {\cal A}(L)$ put $\varepsilon(b):=\varepsilon(A_b)$.  As in the previous section we introduce the space  ${\cal W}^\lambda\subset {\cal A}_\iota(L)\times A^1_\iota(S)$ defined by the inequality $\|\alpha\|_{L^\infty}<\varepsilon(b)$ and its gauge quotient ${\cal Z}^\lambda={\cal W}^\lambda/{\cal G}_\iota$. We obtain:
\begin{pr} The restriction of the instanton moduli problem to the normal neighborhood $\nu({\cal U}^\lambda)$ of the space of $\rho$-twisted, type $\lambda$-reductions ${\cal A}^\lambda_a(E)$  and the restriction of the abelian moduli problem (\ref{twabelian}) to the neighborhood ${\cal W}^\lambda$ of ${\cal A}(L)\times\{0\}$ are equivalent moduli problems.

In particular, the map induced by $\psi$   applies ${\cal M}_{a,\iota}(L)\cap {\cal Z}^\lambda$ isomorphically onto the neighborhood ${\cal M}_a^\ASD(E)\cap {\cal V}^\lambda$ of  the subspace  ${\cal M}_a^\lambda(E)$ of $\rho$-twisted, $\lambda$-reducible instantons, and induces an isomorphism ${\cal M}^{\rm red}_{a,\iota}(L)\simeq {\cal M}_a^\lambda(E)$.

\end{pr}  
\begin{co}  Let $b\in {\cal A}_{\pi_\rho^*(a)}(L)\cap {\cal A}_\iota(L)$. The deformation elliptic complex ${\cal C}_{A_b}$ at the corresponding $\rho$-twisted reduction $A_b$ splits as a direct sum  ${\cal C}_{A_b}={\cal C}_{0,\iota}^+(M_\rho)\oplus {\cal C}_{b\otimes \iota^*(b)^\vee,\iota}^+(S)$, where ${\cal C}_{0,\iota}^+(M_\rho)$ is the  $d^+$-elliptic complex
$$0\map iA^0_\iota(M_\rho)\map iA^1_\iota(M_\rho) \textmap{d^+} iA^2_{+,\iota}(M_\rho)\map  0
$$
of imaginary $\iota$-twisted forms,  and ${\cal C}_{b\otimes \iota^*(b)^\vee,\iota}^+(S)$ is the $d^+$-elliptic complex of $\iota$-twisted $S$-valued forms associated with the connection $b\otimes \iota^*(b)^\vee$.
\end{co}

The indices  of the two  elliptic complexes can be computed easily: The dimension $h^k_\iota$ of the $k$-th harmonic space of ${\cal C}_{0,\iota}^+(M_\rho)$ is  
$$h^0_\iota=0\ ,\  h^1_\iota=b_1(M_\rho)-b_1(M)\ ,\  h^2_\iota=b_+(M_\rho)-b_+(M)\ ,
$$
Taking into account that ${\rm ind}({\cal C}_{0}^+(M_\rho))=2{\rm ind}({\cal C}_{0}^+(M))$, we get
$${\rm ind}({\cal C}_{0,\iota}^+)(M_\rho)=b_+(M_\rho)-b_1(M_\rho)-(b_+(M)-b_1(M))=b_+(M)-b_1(M)+1\ .
$$
For the   complex ${\cal C}_{b\otimes \iota^*(b)^\vee,\iota}^+(S)$, note that the complex ${\cal C}_{b\otimes \iota^*(b)^\vee}^+(S)$ splits as a direct sum
$${\cal C}_{b\otimes \iota^*(b)^\vee,\iota}^+(S)={\cal C}_{b\otimes \iota^*(b)^\vee,0}^+(S)\oplus{\cal C}_{b\otimes \iota^*(b)^\vee,\iota}^+(S)  
$$
where ${\cal C}_{b\otimes \iota^*(b)^\vee,0}^+(S)$ is the $d^+$-complex  of $S$-valued forms $\alpha$ on $M_\rho$ satisfying $\iota^*(\alpha)=\bar\alpha$. Multiplication by $i$ defines a real isomorphism   ${\cal C}_{b\otimes \iota^*(b)^\vee,0}^+(S)\simeq{\cal C}_{b\otimes \iota^*(b)^\vee,\iota}^+(S)$.  Therefore  
$${\rm ind}[{\cal C}_{b\otimes \iota^*(b)^\vee,\iota}^+(S)]={\rm ind}_\C {\cal C}^+(S)= c_1^2(S)+2(b_+(M)-b_1(M)+1)\ .$$
%
 %
%

\subsection{Generic regularity at the reductions}  \label{reg}

The purpose of this section is to prove a strong generic regularity result for reducible instantons. We agree to call {\it regular} any solution $A\in{\cal A}^\ASD_a(E)$ (irreducible or reducible) with $\H^2_A=0$. Our result allows to prove that, under certain cohomological conditions on our data, there exists a {\it connected}, dense, open set of metrics for which no irregular reduction appears in the moduli space.  This will allow us in the next section to introduce Donaldson type invariants for definite manifolds, even in the cases when non-empty reduction loci are present in the moduli space.

Let $M$ be a 4-manifold with $b_+(M)=0$, and $E$ a Hermitian bundle of rank 2 on $M$; put as usual    $D:=\det(E)$, $d:=c_1(D)$, $c=c_2(E)$. Let $\lambda=\{l,d-l\}$ be a topological decomposition of $E$ with $2l\ne d$.  The second cohomology of the deformation elliptic  complex of a $\lambda$-reducible instanton reduces to the second cohomology of the elliptic complex ${\cal C}_{b^{\otimes 2}\otimes a^\vee}^+(S)$, where $S=L^{\otimes 2}\otimes D^\vee$.

Denote  by ${\cal M}et$ (${\cal M}et^r$) the space of smooth (respectively class ${\cal C}^r$) Riemannian metrics on $M$, where $r\gg k$.   Our first result is  a transversality theorem (with respect to variations of the metric $g$) for the complement ${\cal H}_g(S)^*$ of the zero section in the  complex linear space 
$${\cal H}_g(S):=\coprod_{[\sigma]\in{\cal M}^{\ASD_g}(S)} [\H^0({\cal C}^{+_g}_\sigma(S))\oplus \H^2({\cal C}^{+_g}_\sigma(S)]
$$
over the torus  ${\cal M}^{\ASD_g}(S)$.  Unfortunately this result holds only for metrics  having the following property\\ \\
${\bf H}(S):$ {\it The vanishing locus of the $g_0$ harmonic representative of $c_1(S)$ has Hausdorff dimension $\leq 2$.}
\\ \\
This condition is satisfied by any ${\cal C}^\infty$-metric  by a result of B\"ar\cite{B}, and any metric   $g\in {\cal M}et^r_{\geq 2}(c_1(S))$ (see section \ref{admissmet} in the Appendix).\\

The space ${\cal H}_g(S)^*$ can be identified with the ${\cal G}_{x_0}$-quotient  of  the space of solutions of the  system
\begin{equation}\label{irreg}
\left\{\begin{array}{ccc}
F_\sigma^{+_g}&=&0\\
d_\sigma\zeta+d^{*_g}_\sigma\eta&=&0\\
\end{array}\right.
\end{equation}
for triples $(\sigma,\zeta,\eta)\in {\cal A}(S)\times [(A^0(S)\oplus A^2_{+_g}(S))\setminus\{0\}]$.

Regard (\ref{irreg}) as an equation  for  systems $(g,\sigma,\zeta,\eta)$, where $g\in{\cal M}et^r$, $\sigma\in{\cal A}(S)_k$, and $(\zeta,\eta)\in [(A^0(S)_k\oplus A^2_{+_g}(S)_k)\setminus\{0\}]$.  Therefore our configuration space is now 
$${\cal A}_k^*:={\cal A}(S)_k\times[{\cal A}^0(S)_k\oplus  {\cal A}^2_+(S)_k]^*\ ,
$$
where ${\cal A}^0(S)_k$ is the trivial bundle ${\cal M}et^r\times {\cal A}^0(S)_k$ over ${\cal M}et^r$, ${\cal A}^2_+(S)_k$ is the bundle of $S$-valued selfdual forms over ${\cal M}et^r$ (see section \ref{bundles} in the Appendix), and  the symbol $[\ \cdot\ ]^*$ on the right stands  for the complement of the zero-section.
Regarding ${\cal A}_k^*$ as a locally trivial Banach bundle over the Banach manifold ${\cal M}et^r$, and denoting by $p:{\cal A}_k^*\to {\cal M}et^r$ the obvious projection, we see that the left hand terms of the equations  (\ref{irreg}) define  sections $u$, $v$ in the bundles $p^*(i[{\cal A}^2_+]_{k-1})$ and  ${\cal A}_k^*\times A^1(S)_{k-1}$ over ${\cal A}_k^*$.
\begin{thry}\label{trans} The  section $(u,v)$ in the bundle $\left[p^*(i[{\cal A}^2_+]_{k-1})\right]\oplus A^1(S)_{k-1}$ over ${\cal A}_k$ is submersive at any solution $(g_0,\sigma_0,\zeta_0,\eta_0)$ where $g_0$ has the property ${\bf H}(S)$.
\end{thry}
\pf  Note first that, under our assumptions, the connection $\sigma$ cannot admit nontrivial parallel sections, so  applying $d_\sigma^*$ to the second equation, we get $\zeta_0=0$. Use the metric $g_0$ as a background metric to parameterize the manifold ${\cal M}et^r$ and to trivialize the bundles $i[{\cal A}^2_+]_{k-1}$, ${\cal A}^2_+(S)_k$ over this manifold (see section \ref{bundles}). The system (\ref{irreg}) is equivalent  with
\begin{equation}\label{irregnew}
\left\{\begin{array}{ccc}
h^*\{[(h^{-1})^*(F_\sigma)]^+\}&=&0\ \phantom{,}\\
d_\sigma\zeta-*_{g_h}d_\sigma\left[(h^{-1})^*\eta\right]&=&0\ ,\\
\end{array}\right.
\end{equation}
where the upper script $+$ is used now for the selfdual projection with respect to $g_0$. The left hand terms of the equations define a smooth  map 
$$\Gamma({\rm Sym}^+(T_M,g_0))^r\times {\cal A}(S)_k\times [(A^0(S)_k\times A^2_+(S)_k)\setminus\{0\}]\stackrel{(U,V)}{\map} i[A^2_+]_{k-1}\times A^1(S)_{k-1}.
$$
Out task is to prove that the differential of this map at $(\id,\sigma_0,0,\eta_0)$ is surjective. Let $(\alpha,\beta)\in i[A^2_+]_{k-1}\times A^1(S)_{k-1}$ be pair which is $L^2$-orthogonal to the range of this differential.
Using variations of the variables $\zeta$ and $\eta$ (for $h=\id$) we get 
\begin{equation}\label{beta}
d^*_\sigma\beta=0\ ,\  d^+_\sigma\beta=0\ .
\end{equation}

Using the notations of section \ref{bundles} in Appendix, one has
$$\at{\frac{\partial U}{\partial h}}{(\id,\sigma_0,0,\eta_0)}(\chi)=m_-^+(\chi)(F_\sigma)\ ,\ \at{\frac{\partial V}{\partial h}}{(\id,\sigma_0,0,\eta_0)}(\chi)=*d_\sigma(m(\chi)\eta)
$$
Therefore, using variations  $\chi\in A^0({\rm Sym}(T_M,g_0))^r$ of   $h$ with the property 
$$m_-^+(\chi)(F_\sigma)=0\ ,$$
 and noting that $m_+^-(\chi)=[m_-^+(\chi)]^*$ (see section \ref{bundles}), we obtain for any such $\chi$ 
\begin{equation}\label{orth}
0=\langle*d_\sigma(m(\chi)\eta_0),\beta\rangle_{L^2}=-\langle
m(\chi)\eta_0,*d_\sigma\beta\rangle_{L^2}=$$
$$=\langle m_+^-(\chi)\eta_0,d_\sigma\beta\rangle_{L^2} =\langle [m_-^+(\chi)]^*\eta_0,d_\sigma\beta\rangle_{L^2}\ .
\end{equation}
Here we used the fact that the 2-form $d_\sigma\beta$ is ASD. By Remark \ref{iso} in Appendix we see that {\it any} homomorphism $m\in A^0(\Hom(\Lambda^2_-,\Lambda^2_+))^r$ can be written as $m_-^+(\chi)$ for a certain symmetric endomorphism $\chi$. Therefore (\ref{orth}) holds for any such section $m$ for which  $m(F_\sigma)=0$.
 Now regard $\beta$ as an element in   $A^0(\Hom_\R(S^\vee_\R,\Lambda^2_-))_{k-1}$, $\eta_0$ as an element in $A^0(\Hom_\R(S^\vee_\R,\Lambda^2_+))_{k}$,   and denote by $[\eta_0]^*_\R\in A^0(\Hom_\R(\Lambda^2_+,S^\vee_\R))_{k}$ its adjoint with respect to the obvious real inner products. Changing the position of $\eta_0$ in (\ref{orth})  we obtain
$$\left\langle m^*  ,  d_\sigma\beta\circ [\eta_0]^*_\R \right\rangle_{L^2}=0\ ,
$$
for every $m\in A^0(\Hom(\Lambda^2_-,\Lambda^2_+))^r$ for which $m(F_\sigma)=0$. The condition  $m(F_\sigma)=0$ is equivalent to the condition $\im(m^*)\subset F_\sigma^\bot$.
Let $U$ be the complement of the vanishing locus of $F_\sigma$. We conclude that $\resto{d_\sigma\beta\circ [\eta_0]^\vee_\R}{U}$ is $L^2$-orthogonal on the whole space of compactly supported $F_\sigma^\bot$-valued bundle homomorphisms $\resto{\Lambda^2_+}{U}\to \resto{\Lambda^2_-}{U}$.  Therefore  $\resto{d_\sigma\beta\circ [\eta_0]^*_\R}{U}$ takes values in the real line bundle generated by $F_\sigma$. This implies that either there exists a non-empty open subset $V\subset U$ on which  $\eta_0$ has (real) rank at most 1, or  $\resto{d_\sigma\beta}{U}$ takes values in real line bundle generated by $\resto{F_\sigma}{U}$. In the first case we obtain $\eta_0=0$ by Proposition \ref{unkn} in the Appendix. This contradicts  the definition of our configuration space ${\cal A}_k^*$. In the second case one gets $\beta=0$ by Corollary \ref{fub}.
Finally, using variations of $\sigma$ and the assumption  $b_+(M)=0$, we obtain $\alpha=0$.\qed
\\

Denote by ${\cal M}et_{\rm bad}^r(\lambda)$ the subspace of metrics for which there exists a non-regular $\lambda$-reducible instanton. Let $U\subset {\cal M}et^r$ be any open subset of metrics  satisfying the property ${\bf H}(S)$.
 
\begin{thry}  Suppose that   $(2l-d)^2<0$. Then ${\cal M}et_{\rm bad}^r(\lambda)\cap U$ is closed and nowhere dense in $U$ and the natural morphism
$$\pi_i(U\setminus {\cal M}et_{\rm bad}^r(\lambda))\to \pi_i(U)
$$
is  bijective for any $i$ in the range $0\leq i \leq -2(2l-d)^2+b_1(M)-2$ and surjective for $i=-2(2l-d)^2+b_1(M)-1$.
\end{thry}
\pf  By Theorem \ref{trans}   it follows that the section $(u,v)$ is transversal at any solution with metric component in $U$, so the vanishing locus $Z(u,v)\cap p^{-1}(U)$ is a smooth Banach manifold over $U$. The gauge quotient
$$\resto{{\cal H}(S)^*}{U}:=\qmod{Z(u,v)\cap p^{-1}(U)}{{\cal G}_{k+1,x_0}}
$$
will  also be a smooth manifold, and the   natural projection $\resto{{\cal H}(S)^*}{U}\to  {\cal M}et^r$   is Fredholm of real  index $2[c_1(S)^2+(1-b_1(M))]+b_1(M)=2c_1(S)^2-b_1(M)+2$. One has a natural $\C^*$-action on ${\cal H}(S)^*$, and the projection 
$$\qmod{{\cal H}(S)^*}{\C^*}\map  {\cal M}et^r$$
 will be Fredholm of real index $2c_1(S)^2-b_1(M)$. It suffices to apply Lemma  \ref{hgroups} in the Appendix.
\qed

The same arguments can be used to prove regularity at a locus of  twisted reductions (see sections \ref{twred}, \ref{tamp}). However, there is an important detail  here which should be taken into account carefully: in general, for an epimorphism $\pi(M)\to\Z_2$, the condition $b_+(M)=0$ does {\it not} imply $b_+(M_\rho)=0$. When  $b_+(M_\rho)=0$, one   has $h^2_\iota=0$, $h^1_\iota=b_1(M)-1$ (see section \ref{tamp}), and one can obtain generic regularity at the $\rho$-twisted reduction in the same way as for non-twisted reductions, by extending  our proofs  to the twisted case. The main difference is that the ``normal" elliptic complex" ${\cal C}^+_{b\otimes \iota^*(b)^\vee,\iota}$ has no complex structure.    

In the twisted case,  one has a  Hermitian line bundle $S$ on $M_\rho$     which comes with an isomorphism $\iota^*(S)\simeq \bar S$; and for any metric $g$ on $M$ the $\pi_\rho(g)$-harmonic representative of $c_1(S)=l-\iota^*(l)$ is $\rho$-equivariant, so the condition ${\bf H}(S)$ has sense for $g$. Put $\lambda:=\{l,\iota^*(l)\}$ and denote by ${\cal M}et_{\rm bad}^r(\lambda)$ the space of  ${\cal C}^r$-metrics for which there exists a non-regular $\rho$-twisted reducible  instanton of type $\lambda$. The result for the twisted case  is 
\begin{pr}    Suppose that $b_+(M)=b_+(M_\rho)=0$ and let $l\in H^2(M_\rho,\Z)$ a solution of the system (\ref{system-l}) with $(l-\iota^*(l))^2<0$. Put $\lambda:=\{l,\iota^*(l)\}$ and let $U\subset {\cal M}et^r$ an open set of metrics having the property ${\bf H}(S)$.
Then ${\cal M}et_{\rm bad}^r(\lambda)\cap U$ is closed and   the natural morphism
$$\pi_i(U\setminus {\cal M}et_{\rm bad}^r(\lambda))\to \pi_i(U)
$$
is bijective for any $i$ in the range $0\leq i \leq - (l-\iota^*(l))^2+b_1(M)-2$, and surjective for $i=-(l-\iota^*(l))^2+b_1(M)-1$.
\end{pr}
\begin{re} Note that $- (l-\iota^*(l))^2=2[4c_2(E)-c_1(E)^2]$, hence under our assumptions one has $-(l-\iota^*(l))^2\geq 2$.
\end{re}
\begin{co}\label{regularityE} Let $M$ be a 4-manifold $E$ a Hermitian rank 2 bundle on $M$.
\begin{enumerate}
\item Suppose that    
\begin{equation}\label{condition}
b_+(M)=0\hbox{ and } c_1(E)\not\in  2H^2(M,\Z)+\Tors\ .
\end{equation}
 There exists a  connected, dense, open subset ${\cal M}et^r_{\rm good}(E)\subset {\cal M}et^r$ such that, for any $g\in  {\cal M}et^r_{\rm good}(E)$, the reductions in the Uhlenbeck compactification of the moduli space ${\cal M}^\ASD_{a,g}(E)$ of $g$-instantons are all regular. 
\item Suppose that (\ref{condition}) holds and for every epimorphism $\rho:\pi_1(M,x_0)\to \Z_2$ 
\begin{equation}\label{twcondition}
b_+(M_\rho)=0 \hbox{ and } \pi_\rho^*(c_1(E))\not\in 2H^2(M_\rho,\Z)+\Tors  \ .
\end{equation}
 \end{enumerate}
There exists a  connected, dense, open subset ${\cal M}et^r_{\rm vgood}(E)\subset {\cal M}et^r$ such that, for any $g\in  {\cal M}et^r_{\rm vgood}(E)$, the reductions and the twisted reductions in the Uhlenbeck compactification of the moduli space ${\cal M}^\ASD_{a,g}(E)$ of $g$-instantons are all regular. 
\end{co}
\pf 1. Denote by $\Lambda(c)$ the {\it finite} set of unordered pairs $\lambda=\{l,d-l\}$ satisfying $l\cdot(d-l)\leq c$. For $c=c_2(E)$, this set is the set of all topological decompositions of all bundles $E'$ which must be considered in the construction of the Uhlenbeck compactification  of ${\cal M}^\ASD_a(E)$.
Put 
$${\cal M}et^r_{\rm good}(E):={\cal M}et^r_{\rm adm}(c)\setminus \left[\union_{\lambda\in \Lambda(c)} {\cal M}et^r_{\rm bad}(\lambda)\right]$$
(see section \ref{admissmet} in the Appendix).
\\ \\
2. One defines ${\cal M}et^r_{\rm vgood}(E)$ in a similar way by replacing ${\cal M}et^r_{\rm adm}(c)$ with ${\cal M}et^r_{\rm tadm}(c)$ and subtracting  the bad loci associated with all  twisted and non-twisted reductions of bundles $E'$ with $\det(E')=\det(E)$, $c_2(E')\leq c_2(E)$.

\qed
\begin{re} Suppose that $b_+(M)=0$, $b_2(M)>0$, and $b_+(M_\rho)=0$ for any epimorphism $\rho:\pi_1(M,x_0)\to \Z_2$. Let $d$ be an integral lift of $w_2(M)$. Then $d\not\in 2H^2(M,\Z)+\Tors$ and $\pi_\rho^*(d)\not\in 2H^2(M_\rho,\Z)+\Tors(H^2(M_\rho,\Z))$ for every  epimorphism $\rho$, so Corollary \ref{regularityE}.2.  applies for any bundle $E$ with $c_1(E)=d$.
\end{re}
\pf  By Donaldson first theorem, the intersection form on $H^2(M_,\Z)/\Tors$ is trivial over $\Z$. Choosing an orthonormal basis $(e_i)_{1\leq i\leq b_2(M)}$ in this lattice, one obtains $d\cdot e_i\equiv e_i^2\equiv-1$ mod 2, so $d$ cannot be divisible by 2 in $H^2(M,\Z)/\Tors$.

On the other hand, the class $\pi_\rho^*(d)$ is an integral lift of  $\pi_\rho^*(w_2(M))=w_2(M_\rho)$. Since $b_+(M_\rho)=b_+(M)=0$, one gets easily (comparing the signatures and the Euler characteristics of the two manifolds) that $b_1(M_\rho)=2b_1(M)-1$ and $b_2(M_\rho)=2b_2(M)>0$.  Therefore the same argument applies for $M_\rho$, proving that $\pi_\rho^*(d)$ cannot be divisible by 2 in $H^2(M_\rho,\Z)/\Tors$.
\qed
\begin{re} Similar generic  regularity  results can be obtained using abstract perturbations of the ASD-equations  around the reduction loci (see \cite{DK} p. 156). However,  since in our general framework the reductions are not necessarily  isolated points in the moduli space, this method 
is   more complicated than in the classical case. Moreover, for our purposes (see section \ref{appl}) one must check that the perturbed moduli space  still has a natural compactification,  and that the  ``cobordism type" of this compactification  is well defined.
\end{re}

\begin{re}\label{wgood} Combing Corollary \ref{regularityE} with the classical transversality theorem for irreducible instantons \cite{FU}, \cite{DK}, one shows that the  subset   ${\cal M}et^r_{\rm wgood}(E)\subset {\cal M}et^r_{\rm vgood}(E)$ of metrics for which ${\cal M}^\ASD_{a,g}(E)$  contains only regular solutions is dense of the second Baire category. 
This set is also open (but in general non-connected!) when $\Delta(E):=4c_2(E)-c_1(E)^2\leq 3$.
\end{re} 
\pf The condition $\Delta(E)\leq 3$ implies that the projection $\tilde{\cal M}_a^\ASD(E)\to{\cal M}et^r$ of the parameterized instanton moduli space on the space of metrics is proper. The openness of ${\cal M}et^r_{\rm wgood}(E)$ in ${\cal M}et^r_{\rm vgood}(E)$ follows by elliptic semicontinuity.
\qed

 \section{Applications}\label{appl}

\subsection{Geometric properties of instanton moduli spaces on manifolds with $b_+=0$}\label{gp}

We will see that using our regularity results combined with the topological results obtained in sections \ref{ub}, \ref{topology} one can obtain important information about the geometry of the ASD moduli spaces.

The purpose of this section is not to give an exhaustive list of all possible applications of this type, but only the illustrate the method with an explicit example, which came to my attention when I began to work on the classification of class VII surfaces with $b_2=2$ \cite{Te3}.\\

Let $M$ be a 4-manifold with the topological properties considered in       
the examples studied in section \ref{topology}: $H_1(M,\Z)\simeq\Z$, $b_2(M)=b_-(M)=2$.  Consider again a rank 2-Hermitian bundle $E$ on $M$ with $c_1(E)=d=e_1+e_2$ (where $(e_1,e_2)$ is an orthonormal basis of $H^2(M,\Z)\simeq\Z^{\oplus 2}$) and $c_2(E)=0$. Put as in  section \ref{topology} $\lambda_0:=\{0,d\}$, $\lambda_1:=\{e_1,e_2\}$ and note that $\Lambda_E=\{\lambda_0,\lambda_1\}$. Let $L_i$ a Hermitian line bundle of Chern class $d_i$ and $S_i=L_i^{\otimes 2}\otimes D^\vee$.

The expected dimension of the instanton moduli space ${\cal M}^\ASD_a(E)$ is 4 and, since $\Delta(E):=4c_2(E)-c_1(E)^2<4$, this moduli space space is compact. ${\cal M}^\ASD_a(E)$ contains two circles of reductions ${\cal M}^{\lambda_0} _a(E)$ and  ${\cal M}^{\lambda_1} _a(E)$. An interesting application of our results is the following:
\begin{thry} For every Riemannian metric $g$ on $M$ and abelian  connection $a\in{\cal A}(\det(E))$,  the two circles of reductions ${\cal M}^{\lambda_0} _a(E)$, ${\cal M}^{\lambda_1} _a(E)$ belong to the same connected component of the moduli space ${\cal M}^\ASD_a(E)$.
\end{thry}
\pf For a metric $g\in{\cal M}et^r_{wgood}(E)$ the moduli space ${\cal M}^\ASD_a(E)$ contains only regular solutions. Regularity at the reductions implies that the linear spaces 
$${\cal H}_i^1:=\union_{[b]\in T_a(L)} \H^1({\cal C}^+_{b^{\otimes 2}\otimes a^\vee}(S_i))\to 
T_a(L)\simeq {\cal M}^{\lambda_i} _a(E)$$
are rank 2 complex vector bundle of rank 2 (see section \ref{amp}).    ${\cal M}^{\lambda_i} _a(E)$ has a neighborhood  $\nu_i$ which can be identified with the $S^1$-quotient of an $S^1$-invariant  neighborhood of the zero section of ${\cal H}_i^1$. Let $\sigma$ a generator of $H_1(M,\Z)$. By Corollary \ref{dclasses} we see that the restriction of the Donaldson class $\mu(\sigma)$ to the boundary $\partial(\nu_i)$ coincides (up to sign) with the fundamental class  of this 3-manifold. Therefore $\partial(\nu_i)$ cannot be homologically trivial in ${\cal B}^*_a(E)$. This shows that the two boundaries (hence also the corresponding circles) belong  to the same  connected component.\\

To complete the proof for an  arbitrary metric $g$, use the density of the space ${\cal M}et^r_{wgood}(E)$ and note that if the reduction circles ${\cal M}^{\lambda_i} _a(E)$ belonged to different connected components, the same would happen for any metric $g'$ sufficiently close to $g$.
\qed
\\
{\bf Example:} Consider the 4-manifold $M=(S^1\times S^3)\#\bar\P^2\#\bar\P^2$.
This manifold has the differentiable type of a Hopf surface blown up at two points. 
It is convenient to endow $M$ with the complex structure of a {\it minimal} class VII surface with $b_2=2$. Choosing the Gauduchon metric in a convenient way  and using the Kobayashi-Hitchin correspondence to identify instantons with polystable bundles, one obtains (see \cite{Te3}):
$${\cal M}^\ASD_a(E)\simeq S^4\ ,
$$
so (despite the presence of the reductions) the moduli space  gets  an obvious smooth structure on the moduli space. The two reduction  circles ${\cal M}^{\lambda_i} _a(E)$ are smoothly embedded in the sphere.

 \subsection{New Donaldson invariants}
 
 In this section we introduce a new class of Donaldson type invariants, which are defined for definite 4-manifolds.  
 
 \subsubsection{Low energy Donaldson invariants. Casson type invariants. }\label{le}
 
 Let $M$ be a negative definite 4-manifold, and let $(e_1,\dots, e_{b_2(M)})$ be an orthonormal basis in $H^2(M,\Z)/\Tors$.  Let $d\in H^2(M,\Z)$ be a  lift of   $e_1+,\dots,+e_{b_2(M)}$ and denote by $\bar d$ its image in $H^2(X,\Z_2)$. Note that one must have $\bar d=w_2(M)$, when $H^2(M,\Z)$ is torsion free.
  
 Let $E$ be a Hermitian 2-bundle on $M$ with $c_1(E)=d$ and put as usual $D:=\det(E)$. If $\{l,d-l\}$ is a topological decomposition of $E$, then, writing $l=\sum l_i e_i$, with $l_i\in\Z$, one gets  
 $$c_2(E)=\sum l_i(l_i-1)\geq 0\ .
 $$
 Therefore, for $c_2(E)<0$, the bundle $E$ admits no topological decomposition. On the other hand the expected dimension of the Donaldson moduli space  ${\cal M}_a^\ASD(E)$ is
 $$\delta=2(4c_2(E)+b_2(M))+3(b_1(M)-1)\ .$$

If $c_2(E)$ is chosen such that   $\Delta(E)=4c_2(E)+b_2(M)\in\{0,1,2,3\}$  the corresponding moduli space will be a priori compact (i.e. compact independently of the metric).  When $b_2(M)\geq 4$, the corresponding values of $c_2(E)$ are negative. When $b_1(M)\geq 1$ the corresponding expected dimension will be non-negative.
 Therefore:
 \begin{re} Suppose that $b_1(M)\geq 1$ and $b_2(M)\geq 4$, and choose $c_2(E):=-\left[\frac{b_2(M)}{4}\right]$. The corresponding moduli space will be a priori compact, of non-negative expected dimension $2\Delta(E)+3(b_1-1)$ and will contain no reduction.
 \end{re}
 
 In other words, for this special value of $c_2(E)$, one can define very easily Donaldson type invariants by estimating products of classes of the form $\mu(h)$ on the virtual fundamental class of the moduli space \cite{Br}. In this case, one does not really need regular moduli spaces, because  the formalism of  virtual fundamental classes gives directly a well defined homology class in the space ${\cal B}^*_a(E)$.
 
 A very interesting case is when $b_2(M)\geq 4$ is divisible by 4. In this case  this special value of $c_2(E)$ is $-\frac{b_2(M)}{4}$ and the corresponding discriminant $\Delta(E)$ vanishes. Therefore, in this case  ${\cal M}^\ASD(E)$ coincides with the {\it  moduli space  of $PU(2)$-representations of $\pi_1(M,x_0)$ with fixed Stiefel-Whitney class $\bar d$,  modulo $SU(2)$-conjugation.} The invariants associated with such a moduli space should be called four-dimensional Casson type invariants. They should be regarded as an extension of the similar $SU(2)$-invariant  defined for $\Z[\Z]$-homology  $S^1\times S^3$-manifolds (see \cite{RuS}, \cite{FO}) to our new class of homology types.  Note that, because of the absence of reductions in our moduli space, the definition of the invariant in our case is much easier.  These representation spaces can be of course oversized, so it is not clear at all that the corresponding invariants are of  homotopical nature.
   
 Note the following simple, but interesting  vanishing result, which shows that, if non-trivial, this Casson type invariants can be regarded as obstructions to the representability of the basis elements $e_i$ by embedded 2-spheres.
 \begin{re} Suppose that  $b_2(M)$ is divisible by 4 and one of the  basis elements $e_i$ is represented by an embedded sphere. Then the moduli space associated with  the Chern class $c_2(E)=-\frac{b_2(M)}{4}$  is empty. In particular, the corresponding Casson type invariants vanish.
 \end{re}
  If $Z$ is a an embedded surface  representing $e_i$,  one  has $\langle \bar d, [Z]\rangle= \langle  d, [Z]\rangle$ mod 2=1. Let $\rho:\pi_1(M,x_0)\to PU(2)$ a representation of Stiefel-Whitney class $\bar d$. The composition $\pi_1(S,x_0)\to\pi_1(M,x_0)\to PU(2)$ will be a representation with Stiefel-Whitney class $\resto{\bar d}{S}\ne 0$, so $S$ cannot be simply connected.
  \qed
  
  Interestingly, one has:
 \begin{re} 
  There exist  definite negative 4-manifolds, with the property that no element   $e\in H^2(M,\Z)/\Tors$ with $e^2=-1$ can be represented by an embedded sphere.
  \end{re} 
  
  Indeed, it suffices to consider a fake projective plane (see for instance \cite{PY})  with reversed orientation. Since the universal cover of such a 4-manifold is the complex 2-ball, we see that the generator of its homology cannot be represented by an immersed sphere.
\\ 
 
 Consider now the case  
 $$c_2(E)\in \left(-\left[\frac{b_2(M)}{4}\right],-1\right]\ .$$
 In this range, one loses ``a priori compactness", but  has moduli spaces with no reductions in their Uhlenbeck compactifications. In this range, one   uses Donaldson's method    \cite{DK} to define the invariants geometrically:  one uses metrics for which   all strata are regular, constructs distinguished  cycles representing the $\mu$-classes (and which extends to the Uhlenbeck compactification) and defines the invariants by intersecting the moduli space with  systems of such cycles (which can be chosen so that they intersect transversally in the main stratum).

 \subsubsection{Invariants associated with classes in $H^*({\cal B}^*,{\cal V}^*)$. Invariants defined using the cobordism type of the moduli space}

We illustrate these types of invariants in the concrete situation considered in sections \ref{topology} and \ref{gp}: a negative definite 4-manifold with $H_1(M,\Z)\simeq\Z$, $b_2(M)=b_-(M)=2$ endowed with a rank 2-Hermitian bundle $E$ on $M$ with $c_1(E)=d=e_1+e_2$ (where $(e_1,e_2)$ is an orthonormal basis of $H^2(M,\Z)\simeq\Z^{\oplus 2}$) and $c_2(E)=0$.

The moduli space is a priory compact, but it always contains two circle of reductions. We have two ways to define invariants in this situation:\\ \\
1. Use a generic metric in the sense of  \cite{DK}, for which the irreducible part of the moduli space is regular. Regard the (oriented) moduli space ${\cal M}^\ASD_a(E)^*$  of irreducible instantons as a cycle in the relative homology  $H_4({\cal B}^*,{\cal V}^*)$. On the other hand, we have seen in section \ref{topology} that the relative rational cohomology  $H^4({\cal B}^*,{\cal V}^*) \simeq \Q^4$ and that this group fits in a short exact sequence 
\begin{equation}\label{exseq}
0\map \frac{H^3({\cal V}^*_0)\oplus  H^3({\cal V}^*_1)}{\langle(\delta_{L_0}(\sigma_0) h_{u_0},\delta_{L_1}(\sigma_0) h_{u_1} ) \rangle}\map H^4({\cal B}^*,{\cal V}^*) \map \ \ \ \ \ \ \  \ \ \ \ \ \ \  \ \ \ \ \ \ \   $$
$$\ \ \ \ \ \ \    \ \ \ \ \ \ \  \ \ \ \ \ \ \   \map\langle \mu(*)+\mu(f_1)^2, \mu(*)+ \mu(f_2)^2,\ \mu(\lambda_0)\mu(\sigma_0)\rangle_\Q\to 0\ .
\end{equation}

Evaluating classes in  $H^4({\cal B}^*,{\cal V}^*)$ on the relative homology class defined by ${\cal M}^\ASD_a(E)^*$ one gets well-defined invariants. Note however that, since the exact sequence (\ref{exseq}) does not split canonically, one cannot parameterize this set of invariants in an obvious way.
\\ \\
2. For a metric $g\in {\cal M}et^r_{\rm wgood}(E)$ all   solutions (including the reductions) in the moduli space are regular (see Remark \ref{wgood}).    The main observation here is that each reduction circle ${\cal M}^{\lambda_i}_a(E)$  has a  neighborhood  isomorphic to the $S^1$-quotients of a neighborhood of  the zero section in a rank 2 complex bundle ${\cal H}^1_i$ over  ${\cal M}^{\lambda_i}_a(E)$ (see section \ref{gp}). But such a quotient is a locally trivial $K$-bundle over a circle, where $K$ is the  cone over $\P^1\simeq S^2$, so it has a natural manifold structure.
Therefore, for $g\in {\cal M}et^r_{\rm wgood}(E)$, ${\cal M}^\ASD_a(E)$ is a compact 4-manifold, which can be oriented as in classical Donaldson theory (see \cite{DK} p. 283).

The signature of this 4-manifold will be an invariant $\theta(M)$ of the base 4-manifold (endowed with the usual orientation data).  Indeed, the main point here is that the set ${\cal M}et^r_{\rm vgood}(E)$ is connected. For two choices $g_0$, $g_1\in {\cal M}et^r_{\rm wgood}(E)$, consider a path $\gamma:[0,1]\to {\cal M}et^r_{\rm vgood}(E)$ connecting these metrics. A generic deformation (with fixed ends) of $\gamma$ will define a cobordism between the moduli spaces associated with $g_i$. Note that  {\it the cobordism  constructed in this way is always trivial around the reductions}. 
\begin{re} In the example given in section \ref{gp}, one has ${\cal M}^\ASD_a(E)\simeq S^4$, so the $\theta$-invariant vanishes. Using gluing theory for instantons, one can prove more generally that this invariant also vanish for a manifold of the form $M=N\#\bar\P^2$, where $N$ is a 4-manifold with $H_1(N,\Z)\simeq \Z$ and $b_2(N)=b_2^-(N)=1$.
\end{re}

This suggests that, if non-trivial, this invariant  can also be regarded as an obstruction to the representability of the elements of the orthonormal basis $\{e_1,e_2\}$ by embedded spheres.

\section{Appendix}

\subsection{Metric-dependent spaces of selfdual and anti-selfdual forms} \label{bundles}

Let $M$ a compact oriented  connected 4-manifold and $E$ a vector bundle. For every Riemannian metric $g$ on $M$ one has two associated spaces of $E$-valued (anti)selfdual forms $A^2_{\pm_g}(E)$. It is convenient to complete the space ${\cal M}et$ of metrics with respect to the ${\cal C}^r$ topology and the spaces $A^2_{\pm_g}(E)$ with respect to a Sobolev norm $L^2_k$ (where $r\gg k$). In this way one gets Banach vector bundles $[{\cal A}^2_\pm(E)]_k$  on the Banach manifold ${\cal M}et^r$ of ${\cal C}^r$-metrics. 

One can trivialize globally these bundles in the following way. Fix a ${\cal C}^r$-metric $g_0$.  The space ${\cal M}et^r$ can be identified with the space of positive  $g_0$-symmetric automorphisms of  the tangent bundle $T_M$ via the diffeomorphism $h\mapsto g_h:=h^*(g_0)$.  We get homeomorphisms
$$\Gamma({\rm Sym}^+(T_M,g_0))^r\times   [A^2_{\pm_{g_0}}(E)]_k\textmap{\simeq} [{\cal A}^2_\pm(E)]_k 
$$
given by $(h,\eta)\mapsto (g_h,h^*(\eta))$.  It is important to notice that   homeomorphisms associated with different metrics $g_0$ are pairwise differentiable compatible. Therefore one can use these homeomorphisms to define   structures of Banach manifolds  on  the total spaces $[{\cal A}^2_\pm(E)]_k$.

A positive symmetric automorphism $h\in \Gamma({\rm Sym}^+(T_M,g_0))^r$  defines a  class ${\cal C}^r$ positive $g_0$-symmetric automorphism $\Lambda^2h$ of the bundle $\Lambda^2_M$, given by $\lambda\mapsto \chi^*(\lambda)$. Using the $g_0$-orthogonal decomposition  $\Lambda^2_M=\Lambda^2_{+_{g_0}}\oplus\Lambda^2_{-_{g_0}}$, the automorphism $\Lambda^2h$ can be written as
$$\Lambda^2h=\left(
\begin{array}{cc}
\mu_+(h)&\mu^+_-(h)\\
\mu^-_+(h)&\mu_-(h)
\end{array}\right)\ ,
$$
where $\mu_\pm(h)$ is a symmetric endomorphisms of $\Lambda^2_{\pm_{g_0}}$ and 
$\mu^\pm_{\mp}(h):\Lambda^2_{\mp_{g_0}}\to \Lambda^2_{\pm_{g_0}}$ have the property $\mu^+_-(h)^*=\mu^-_+(h)$. The tangent space of $\Gamma({\rm Sym}^+(T_M,g_0))^r$ at $\id$ can be identified with $A^0({\rm Sym}(T_M,g_0))^r$. For a symmetric endomorphism $\chi$ we put 
$$m^\pm_\mp(\chi):=\at{\frac{\partial m^\pm_\mp(h)}{\partial h}}{\id}(\chi)\ .$$
 One also has $m^+_-(\chi)^*=m^-_+(\chi)$.
 \begin{re}\label{iso} The map $\chi\mapsto m^\pm_\mp(\chi)$ defines an isomorphism  
 $$A^0({\rm Sym}_0(T_M,g_0))^r\simeq A^0(\Hom(\Lambda^2_\mp,\Lambda^2_\pm))^r\ , $$
  where ${\rm Sym}_0(T_M,g_0)$ denotes the bundle of trace free $g_0$-symmetric endomorphisms of $T_M$.
 \end{re}
 
 The space $A^0({\rm Sym}(T_M,g_0))^r$ can be regarded as the space of infinitesimal variations of ${\cal C}^r$-metrics, whereas $A^0({\rm Sym}_0(T_M,g_0))^r$ is the space of infinitesimal variations of ${\cal C}^r$-conformal structures.
 
 \subsection{Admissible metrics} \label{admissmet}

Consider an oriented  compact connected manifold $M$ and a cohomology class $l\in H^p_{\rm DR}(M)$, where $1\leq p\leq n-1$. For a Riemannian metric $g$ on $M$ we will denote by $l_g$ the unique  $g$-harmonic representative of $l$. 

Let $M$ a 4-manifold. The family of vector spaces 
$$(\Lambda^2_{\pm_g,x})_{g\in {\cal M}et^r, x\in M}$$
 defines a bundle  $\Omega^\pm$ on the product space ${\cal M}et^r\times M$.
Suppose  now that $b_+(M)=0$. In this case any harmonic 2-form is ASD. For a de Rham 2-cohomology class $l$, the assignment ${\cal M}et^r\times M\ni (g,x)\mapsto l_g(x)$ defines a universal  ${\cal C}^{r,\epsilon}$ section $\lambda$ in the the bundle $\Omega^-$ over ${\cal M}et^r\times M$.

\begin{lm}\label{transv}   Suppose that $b_+(M)=0$ and $l\in H^2_{\rm DR}(M,\R)\setminus\{0\}$.   The associated universal section $\lambda$ is  submersive at every vanishing point.
\end{lm}
\pf Let $(g_0,x_0)$ be a a vanishing point of $\lambda$. We use the metric $g_0$ to simultaneously parameterize the space of metrics (using positive $g_0$-symmetric endomorphisms) and to trivialize the bundles  $\Omega^{\pm}$.
 In this way, our section gives rise to a map 
 $$\Gamma({\rm Sym}^+(T_M,g_0))^r\times M \ni (h,x)\mapsto  \omega_h(x)\in \Lambda^2_{-,x} \ .$$
where $\omega_h:=[\mu(h)^{-1}]^*(l_{g_h})\in A^2_-(M)$. Put $\omega_0:=\omega_{\id}=l_{g_0}$. Writing $l_{g_h}=\omega_0+d\alpha_h$ (with $\alpha_h\bot_{L^2} Z^1$), one gets  for the derivative $\dot\alpha(\chi):=\at{\frac{d\alpha}{dh}}{\id}(\chi)$ the equation
$$p^+\left[-m(\chi)^*(\omega_0)+d\dot\alpha(\chi)\right]=0\ ,
$$
which gives the solution
$$\dot\alpha(\chi)=G^+(m(\chi)^+_-(\omega_0))\ ,$$
where $G^+$ is the $[Z^1]^\bot$-valued inverse of the operators $d^+$.
Therefore 
$$\at{\frac{d\omega_h}{d h}}{\id}(\chi)=dG^+(m^+_-(\chi)(\omega_0))-m(\chi)(\omega_0)= d^-G^+(m^+_-(\chi)(\omega_0))-m^-_-(\chi)(\omega_0)\ . 
$$
Therefore
$$\at{\frac{\partial \omega_h(x)}{\partial h}}{(\id,x_0)}(\chi)=d^-G^+(m^+_-(\chi)(\omega_0))(x_0)\ ,
$$
because $\omega_0(x_0)=0$.
Suppose now that $v\in  [\Lambda^-_{x_0}]$ is orthogonal  on  the range of this partial derivative.  The vector $v$ defines a Dirac type $\Lambda^-$-valued distribution $v_{x_0}$. Since any element of $m\in A^0(\Hom(\Lambda^2_-,\Lambda^2_+))$ has the form $m^+_-(\chi)$, we obtain that for any such $m$
\begin{equation}\label{allm}
0=\langle v,d^-G^+(m(\omega_0))(x_0) \rangle=\langle v_{x_0}, d^-G^+(m(\omega_0))\rangle=$$
$$=\langle  (d^-)^*(v_{x_0}), G^+(m(\omega_0)) \rangle
\end{equation}

The assignment $\psi\stackrel{\delta}{\mapsto} \langle (d^-)^*(v_{x_0}), G^+(\psi)\rangle$
is a  $\Lambda^2_+$-valued distribution, and the identity  (\ref{allm}) implies that this distribution vanishes on $M\setminus Z(\omega_0)$.  This distribution has the following important properties:

\begin{enumerate}
\item One has $(d^+)^*(\delta)=(d^-)^*(v_{x_0})$.
\item $\delta$ is supported in $x_0$.
\end{enumerate}
1. Indeed, to prove  the first claim, consider   a $\Lambda^1$-valued test function $\alpha$, and compute
$$\langle (d^+)^*(\delta),\alpha\rangle =\langle \delta, d^+\alpha\rangle= \langle (d^-)^*(v_{x_0}), G^+d^+\alpha\rangle=\langle(d^-)^*(v_{x_0}), {\rm pr}_{{Z^1}^\bot}(\alpha)\rangle$$
$$=\langle v_{x_0},  d^- ({\rm pr}_{{Z^1}^\bot}(\alpha))\rangle=\langle(d^-)^*(v_{x_0}), d^-\alpha\rangle=\langle (d^-)^*(v_{x_0}), \alpha\rangle\ .$$
2. Using the fact that $(d^+)^*$ is overdetermined elliptic on $A^2_+(M)$,  it follows  from the first statement  hat $\resto{v_{x_0}}{M\setminus\{x_0\}}$ is smooth. On   the other hand $\delta$   vanishes on the  dense open set  $M\setminus Z(\omega_0)$. Therefore $\delta$ must vanish everywhere on $M\setminus\{x_0\}$.
\\

Since $\delta$ is SD and $v_{x_0}$ is ASD, the first property implies $d(\delta+v_{x_0})=0$. The statement follows now  from the following Lemma.
\qed
\begin{lm} Let $u$ be a $\Lambda^2$-valued distribution on $\R^4$ which is supported at the origin and is closed. If $u^-$ has order  0, then $u=0$.
\end{lm}
\pf 
If $du=0$ and $u^-$ has order 0, then $u^+$ must also have order 0, because, if not, taking the  sum $u^+_k$ of all terms of highest order $k\geq 1$  in the decomposition of $u^+$ as sum of partial derivatives  of Dirac type distributions, one would have $du^+_k=0$. Therefore $u^+_k$ would be singular harmonic self-dual form, which is impossible.  

Therefore, $u$ is an order zero Dirac type distribution; in other words, there exists $\theta_0\in [\Lambda^2_0]$ such that $\langle u,\varphi\rangle=\theta_0(\varphi(0))$ for any test 2-form $\varphi$. Writing $\theta_0=\sum_{i<j} a_{ij} e^i\wedge e^j$, we get  
$$du=\sum_{i<j<k} (a_{jk}\frac{\partial}{\partial x_i}-a_{ik}\frac{\partial}{\partial x_j}+a_{ij}\frac{\partial}{\partial x_k}) e^i\wedge e^j\wedge e^k\ ,
$$
and the relation $du=0$ implies obviously $a_{ij}=0$ for all $i<j$.
\qed

\begin{dt}  Let $M$ be a 4-manifold with $b_+(M)=0$ and $l\in H^2_{\rm DR}(M)\setminus\{0\}$.  A metric $g\in {\cal M}et^r$ will be called 
\begin{enumerate}
\item strictly $l$-admissible, if it is submersive (transversal to the zero section) at any vanishing point. 
\item  $l$-admissible, if the intrinsic derivative of the section $l_g\in \Gamma(\Lambda^2_{-,g})$ at any vanishing point has rank at least 2.
\end{enumerate}
\end{dt}

Denote by ${\cal M}et^r_3(l)$ and   ${\cal M}et^r_{\geq 2}(l)$ the space of (strictly) admissible ${\cal C}^r$-metrics.

\begin{pr}\label{propadmiss} \begin{enumerate}
\item Suppose that the regularity class $r$ is sufficiently large. The space ${\cal M}et^r_3$ of $l$-strictly admissible  metrics is  open and dense   in ${\cal M}et^r$.
\item The space 
${\cal M}et^r_{\geq 2}$ of $l$-admissible  metrics is  open, dense and path connected  in ${\cal M}et^r$.
\end{enumerate}
\end{pr}
\pf The openness of the two sets is obvious taking into account the compactness of the manifold, the continuity of the map $g\mapsto l_g$ with respect to the ${\cal C}^1$ topology on the space of sections, and the fact that, in general, for any bundle $E$ and fixed base point $x$, the  condition 
$$s(x)=0\ ,\ \rk(D_x(s))\leq k$$
is closed with respect to the ${\cal C}^1$-topology on the space of sections  in $E$.

To prove the density of the set of admissible metrics, note that the vanishing locus  $Z(\lambda)$ of the universal section $\lambda$ is a smooth codimension 3-submanifold of  ${\cal M}et^r\times M$, and the natural map $q:Z(\lambda)\to {\cal M}et^r$ is Fredholm of index 1. ${\cal M}et^r_3(l)$ is just the set of  regular values of $q$.

Therefore the larger set ${\cal M}et^r_{\geq 2}(l)$ will also be dense. In order to prove the fact that this set is also path connected, consider --  for  two given metrics $g_0$, $g_1\in {\cal M}et^r_{\geq 2}(l)$ -- a differentiable path $\gamma:[0,1]\to {\cal M}et^r$ joining them. 

Consider that map
$$G:\Gamma(SL(T_X))^r\times [0,1] \map {\cal M}et^r
$$
defined by $G(h,t)=h^*(\gamma(t))$.  Since the partial derivative $\frac{\partial G}{\partial h}$ alone  is  surjective at any point, we  conclude that the map $G$ is transversal to $q$, so the fibred product  
$P:=\{(h,t,z)| G(h,t)=q(z)\}$
 is  a smooth submanifold of $\Gamma(SL(T_X))^r\times [0,1]\times  Z$.  The natural map $P \to \Gamma(SL(T_X))^r$ is proper and Fredholm of index 2. By Sard-Smale theorem, in any neighborhood of $\id$ one can find a regular value $h_0$ of this map. Therefore, for such $h_0$  
 $$P_0:=\{(t,z)|\  G(h_0,t,z)=q(z)\}
 $$
is a smooth 2-dimensional submanifold of $[0,1]\times   Z$. Put $g'_t:=G(h_0,t)$. The fibre over $t\in[0,1]$ is contained in $\{t\}\times \{g'_t\}\times M$ and is identified with the  vanishing locus $Z(l_{g'_t})$ of the section $l_{g'_t}$ under the projection on the third factor.
 
 We claim that the path $G(h_0,\cdot)$ takes values in ${\cal M}et^r_{\geq 2}(l)$.  
  Indeed, the Zariski tangent space of  $Z(l_{g'_t})$ at a   point $x$ is just the intersection of the tangent space of $P_0$ at $(t,g'_t,x)$ with the tangent space of the fibre $\{t\}\times \{g'_t\}\times M\subset P_0$ over $t$. But the tangent space of $P_0$ at any point is 2-dimensional, so the Zariski tangent space of  $Z(l_{g'_t})$ at   $x$ has dimension 1 or 2.
  
  In order to complete the proof, is suffices to join $g_0$ to $g_0'$ and $g_1$ to $g_1'$ with paths in  ${\cal M}et^r_{\geq 2}(l)$. By the openness property of this space, it follows easily that ( if  $h_0$ is sufficiently close to $\id$)  the metrics $G(\id+s(h_0-\id,0)$, $G(\id+s(h_0-\id,0)$, $s\in[0,1]$ remain in    ${\cal M}et^r_{\geq 2}(l)$ for every $s\in [0,1]$.\qed

  Note that we see no reason why the space  ${\cal M}et^r_{3}(l)$ should be  connected. Indeed, suppose that for two metrics $g_0$, $g_1\in {\cal M}et^r_{3}(l)$, the corresponding vanishing loci (which are finite unions of  pairwise disjoint embedded circles)   have different number of connected components. Then there is certainly no way to join the two metrics by a path in ${\cal M}et^r_{3}(l)$
  
 \begin{re} 
If $g\in {\cal M}et^r_{\geq 2}(l)$, then any point $x\in Z(l_g)$ has a neighborhood $U_x\subset M$ such that $Z(l_g)\cap U$ is contained in a closed submanifold $N_x$ of $U_x$ of dimension 1 or 2. In particular the vanishing locus $Z(l_g)$ of the $g$-harmonic ASD form $l_g$ has Hausdorff dimension $\leq 2$.
\end{re}
\pf Indeed, consider a smooth map $\R^n\supset V\textmap{f}\R^m$ whose rank at  $0\in V$ is $k$ and let $E:=\im(d_0(f))$. Then the composition $p_E\circ f$ is a submersion at $0$, so its restriction on a sufficiently small neighborhood $U$ of $0$ will   be a submersion. Note that $Z(f)\subset Z(p_E\circ f)$ and that $Z(p_E\circ f)\cap U=Z(\resto{(p_E\circ f)}{U})$ is a codimension $k$ submanifold of $U$.
\qed

The results above can be extended for twisted de Rham cohomology classes: Let $\rho:\pi_1(M)\to \Z_2$ an  epimorphism, and suppose that $b_+(M_\rho)=0$ where $\pi_\rho:M_\rho\to M$ is the corresponding double cover of $M$. Let   $l\in H^2_{\rm DR}(M_\rho)\setminus\{0\}$ be a de Rham cohomology class with the property that $\iota^*(l)=-l$, where $\iota$ stands for the tautological involution of $M_\rho$.   Then, for every metric $g$ on $M$, the $\pi_\rho^*(g)$-harmonic representative $l_g$ of $l$ is an ASD form on $M_\rho$ satisfying the identity $\iota^*(l_g)=-l_g$. In other words, $l_g$ is a $\rho$-twisted ASD form on $M$. At any point $x\in M$, $l_g$ is defined up to sign. Therefore, the vanishing locus $Z(l_g)\subset M$ and the rank of the intrinsic derivative at a vanishing point are well-defined objects.  In particular one can associate to $l$ the sets  of metrics ${\cal M}et^r_{3}(l)$, ${\cal M}et^r_{\geq 2}(l)$ as in the non-twisted case, and these sets  also have the properties stated in Proposition  \ref{propadmiss}.
\begin{pr}\label{admiss(c)} Let $M$ be a 4-manifold and $d\in H^2(M,\Z)$ %
\begin{enumerate}
\item Suppose that    
\begin{equation}\label{condition-app}
b_+(M)=0\hbox{ and } d\not\in  2H^2(M,\Z)+\Tors\ .
\end{equation}
 Then, for every $c\in H^4(M,\Z)$ the set 
 $${\cal M}et^r_{\rm adm}(c):=\bigcap_{\matrix{\scriptstyle l\in H^2(M,\Z)\cr\scriptstyle l\cdot(d-l)\leq c}}   {\cal M}et^r_{\geq 2}(2l-d)$$ 
 is open, dense and connected in  ${\cal M}et^r$. 
\item Suppose that (\ref{condition-app}) holds and, for every epimorphism $\rho:\pi_1(M,x_0)\to \Z_2$ 
\begin{equation}\label{twcondition-app}
b_+(M_\rho)=0 \hbox{ and } \pi_\rho^*(d)\not\in 2H^2(M_\rho,\Z)+\Tors  \ .
\end{equation}
 \end{enumerate}
There the set
$${\cal M}et^r_{\rm tadm}(c):={\cal M}et^r_{\rm adm}(c)\bigcap\left[ \bigcap_{\matrix{\scriptstyle  \pi_1(M,x_0)\stackrel{\rho}{\twoheadrightarrow} \Z_2,\    l\in H^2(M_\rho,\Z)\cr\scriptstyle l+\iota^*(l)=\pi_\rho^*(d),\ l\cdot \iota^*(l)\leq \pi_\rho*(c)}}  {\cal M}et^r_{\geq 2}(l-\iota^*(l))\right]$$ 
 is open, dense and connected in  ${\cal M}et^r$. 
\end{pr}
\pf  The important point here  is that, since the intersection forms of $M$ is negative definite, the set of terms in the first intersection is finite. For the second  intersection note first that $\rho$ varies in a finite set  (isomorphic to the set of index 2 subgroups of $H_1(M,\Z)$) and, under our assumption, for any fixed $\rho$, there are only finitely many possibilities for $l$. To complete the proofs it suffices to note that the two properties in Proposition \ref{propadmiss} have been obtained by applying the Sard-Smale theorem to certain proper Fredholm maps, and using the fact that the set of regular values of such a map is open and dense. A finite intersection of such sets will also be open and dense.
\qed

\subsection{Analytic results}

We begin with the following easy result concerning the  range of a proper Fredholm map $f:V\to W$ of negative index $j$. One expects this range to be a subset of    ``codimension" at least $-j$, so  the inclusion $W\setminus\im(f)\to  W$ should induce isomorphism between homotopy groups of degree $i\leq -j-2$. While the statement about the codimension is difficult to formalize in the infinite dimensional framework (and is not necessary for our purposes), the isomorphism of  homotopy groups of small degree is a simple consequence of the Sard-Smale theorem. 
\begin{lm} \label{hgroups}
 Let $V$, $W$ be separable Banach manifolds and $f:V\to W$ a proper, real analytic, Fredholm map of negative index $j\leq -2$. Then
\begin{enumerate}
\item The range  $f(V)$ of $f$ is closed and nowhere dense,
\item The natural map  $\pi_i(W\setminus f(V))\to \pi_i(W)$ is an isomorphism for $0\leq i\leq -j-2$ and is surjective for $i=-j-1$.
 %
\end{enumerate}
\end{lm}
\pf  Use the Sard-Smale theorem and  standard transversality arguments   (see \cite{DK}) to approximate maps $S^i\to W$ and homotopies between such maps by maps (respectively homotopies) with values in $W\setminus f(V)$. It suffices to note that, in general, for $i\leq -j-1$, a map $\phi: N\to W$ from a $i$-dimensional manifold $N$ is transversal to $f$ if and only if $\phi(N)\cap f(V)=\emptyset$.
\qed

The following proposition plays an important role in the proof of our regularity results.  This will allow us to refine the Freed-Uhlenbeck theorems \cite{FU} and to estimate the codimension of the spaces of bad metrics.

\begin{pr}\label{char} Let $M$ be an $n$-dimensional compact manifold, $E$ a  real  rang $r$ vector bundle and $\nabla$ a linear connection on $E$. Let $U\subset M$ be an open set whose complement $\Sigma$ has   Hausdorff dimension  $d\leq n-2$. Let $\alpha\in A^1(E)$ and $\varphi\in\Gamma (U,E)$ such that $\nabla\varphi=\resto{\alpha}{U}$. Then $\varphi$ extends smoothly to a section  $\psi\in \Gamma(M,E)$ satisfying $\nabla\psi=\alpha$.
\end{pr}
The idea of the proof is to use the (classical) methods of characteristics to solve the first order equation $\nabla f=\alpha$: we integrate  the family of ordinary equations obtained by restricting our equation to an $(n-1)$-dimensional family of embedded paths which define a local foliation.  We choose the paths such that all starting points belong to $U$, and we use the values of $\varphi$ at these points as initial conditions. The condition on the Hausdorff dimension of $\Sigma$ implies  that a dense family of paths do not meet the set $\Sigma$ (where $\varphi$ is not defined), so on these paths the sections $f$ and $\varphi$ coincide.  In this way one checks that the section $f$ obtained by the method of characteristics agrees with $\varphi$ on the intersections of their domains, providing a proper extension of $\varphi$. \\ \\
\pf  
 For every smooth path $\gamma:(-1,1)\to M$, consider the connection $\gamma^*(\nabla)$ on the bundle $\gamma^*(E)$ on $(-1,1)$ and the affine ordinary differential equation
\begin{equation}\label{ord}
\gamma^*(\nabla)g=\gamma^*(\alpha)
\end{equation}
for sections $g$ in the bundle $\gamma^*(E)$. Using the general theory of ordinary differential equations, one gets, for every $e\in E_{\gamma(0)}$, a unique solution
$$g_{\gamma,e}\in \Gamma((-1,1),\gamma^*(E))
$$
of the equation (\ref{ord}) satisfying the initial condition $g_{\gamma,e}(0)=e$. This solution depends differentiably on the pair $(\gamma,e)$, where $e\in E_{\gamma(0)}$. Our hypothesis  $\nabla\varphi=\resto{\alpha}{U}$ implies $\varphi\circ\gamma=g_{\gamma,\varphi(\gamma(0))}$
for every smooth path $\gamma:(-1,1)\to U$.  

Let $\psi\in \Gamma(V,E)$ be a maximal element of the ordered set  of extensions of $\varphi$ defined on open subsets of $M$. The existence of such a maximal element follows by the Zorn lemma. Since  $U$ is dense in $M$ (so also in $V$) one has
\begin{equation}\label{eqpsi}
\nabla\psi=\resto{\alpha}{V}
\end{equation}
which implies 
\begin{equation}\label{solution}
\psi\circ\gamma=g_{\gamma,\psi(\gamma(0))}
\end{equation}
for every smooth path $\gamma:(-1,1)\to V$. The complement $\Sigma'=M\setminus V $ will also have Hausdorff dimension $\leq n-2$.

We claim that $V=M$. Suppose not, and let $x_0\in M\setminus V$.  Let 
$$h:B^n(0,1)\times (-1,1)\textmap{\simeq} W\subset M
$$
a local parametrization of $M$ with the properties 
\begin{enumerate}
\item $h(B^n(0,1)\times\{0\})\subset V$, $u(0,\frac{1}{2})=x_0$,
\item $h$ is the restriction of a local parametrization  
$$\tilde h:  B^{n-1}(0,1+\varepsilon)\times (-1-\varepsilon,1+\varepsilon)\map  M\ .$$
\end{enumerate} 
The second property implies that $h$ and $h^{-1}$ are Lipschitz with respect to any Riemannian metric on $M$.  The idea is to extend $\psi$ on $V\cup W$ using the solutions $g_{h_x,\psi(x,0)}$, $x\in B^{n-1}$, where $h_x$ denotes the path $t\mapsto h(x,t)$.

So put
$$f(x,t)=g_{h_x,\psi(x,0)}(t)\ ,\ \eta:=f\circ h^{-1}\ .
$$
We claim that $\resto{\eta}{W\cap V}=\resto{\psi}{W\cap V}$, which will complete the proof, because this would yield a proper extension of $\psi$, contradicting its maximality.

Our claim is equivalent to

$$\resto{f}{h^{-1}(V)}=\resto{\psi\circ h }{h^{-1}(V)}\ .$$
The two functions coincide on $\left[B(0,1)\setminus {\rm pr}_1(h^{-1}(\Sigma'))\right]\times (-1,1)$, because, for any $x\in \left[B(0,1)\setminus {\rm pr}_1(h^{-1}(\Sigma'))\right]$ the corresponding path $h_x$ is entirely contained in  $V$ so, for such $x$, both sections
$$t\mapsto f(x,t)\ ,\ t\mapsto \psi(u(x,t))
$$
 coincide with $g_{h_x,\psi(x,0)}$, by (\ref{solution}) and the definition of $f$.
It suffices to notice that, by  our hypothesis, the set $h^{-1}(\Sigma')$ is of Hausdorff dimension at most $n-2$, so its projection on the $(n-1)$ dimensional ball $B(0,1)$ is also of Hausdorff dimension at most $n-2$. Therefore, this projection cannot contain any non-empty open set,  so its complement in  $B(0,1)$ is dense, so  $\left[B(0,1)\setminus {\rm pr}_1(h^{-1}(\Sigma'))\right]\times (-1,1)]$ is dense in $B^n(0,1)\times (-1,1)$, so the two functions coincide everywhere.
\qed

 The following corollary shows that the statement of  Lemma 4.16 in \cite{FU} is true as soon as the vanishing of the ASD curvature has Hausdorff dimension $\leq 2$. In particular this statement is true for ${\cal C}^\infty$-metrics by the results of \cite{B} and for admissible metrics.

\begin{co}\label {fub} Let $(M,g)$ be  compact oriented Riemannian 4-manifold endowed with a ${\cal C}^r$-metric $g$, and $S$ a Hermitian line bundle on $M$  endowed with a non-flat  ASD connection $\sigma$. Suppose that the vanishing of the curvature has Hausdorff dimension $\leq 2$. Let $\beta\in A^1(S)$ such that 
\begin{enumerate}
\item $d^*_\sigma\beta=0$, $d^+_\sigma\beta=0$,
\item The anti-selfdual $S$-valued form $d_\sigma\beta$ is a tensor multiple of $F_\sigma$ at any point $x\in M$ for which $F_{\sigma,x}\ne 0$.  
\end{enumerate}
Then $\beta=0$.
\end{co}
\pf  Let $U$ be the complement of the vanishing locus $\Sigma$ of $F_\sigma$. By assumption we can write $\resto{d_\sigma\beta}{U}=F_\sigma\otimes\zeta$, for a section $\zeta\in\Gamma(U,S)$. We get
$$[\resto{F_\sigma}{U}]\wedge(\resto{\beta}{U}-d_\sigma \zeta)=d_\sigma d_\sigma(\resto{\beta}{U}-d_\sigma \zeta)=d_\sigma (d_\sigma \resto{\beta}{U}-F_\sigma \otimes\zeta)=0\ ,
$$
hence $\resto{\beta}{U}-d_\sigma \zeta=0$, since the wedge product with a non-trivial ASD form is invertible on 1-forms. By Proposition \ref{char} and the assumption on $Z(F_\sigma)$, the section $\zeta$ extends smoothly to a section $\xi$ on $M$ satisfying $d_\sigma\xi=\beta$. Since $d_\sigma^*\beta=0$, we get immediately $\beta=0$.
\qed

\begin{pr}\label{unkn} Let $S$ be a Hermitian line bundle on an oriented Riemannian 4-manifold  $(M,g)$ and $\sigma$ a  Hermitian connection on $S$.  Let $\eta\in A^2_+(S)$ with $d_\sigma\eta=0$. Suppose that on an open set $U\subset M$, the form $\eta$ (regarded as section in $\Lambda^2_+\otimes S$) has real rank 1. Then $\resto{F_\sigma}{U}=0$.
\end{pr}
\pf Supposing that $U$ is simply connected, we can write $\eta=\omega\otimes \zeta$, where $\omega\in A^2_+(U)$   is a real selfdual form, and $\zeta\in \Gamma(U,S)$. By assumption, both $\omega$ and $\zeta$ are nowhere vanishing  on $U$. Since $d_\sigma\eta=0$, we get
$$d\omega\otimes \zeta+\omega\wedge d_\sigma\zeta=0\ ,
$$
hence $\omega\wedge(\theta\otimes\zeta+ d_\sigma\zeta)=0$, where $\theta$ is the real form on $U$ defined by $d\omega=\omega\wedge\theta$. We get
$$d_\sigma\zeta=-\theta\otimes\zeta\ ,\ F_\sigma\otimes\zeta=-d\theta\otimes\zeta+\theta\wedge d_\sigma\zeta=-d\theta\otimes\zeta-(\theta\wedge  \theta)\otimes\zeta=-d\theta\otimes\zeta\ .
$$
This yields $F_\sigma=d\theta$, in which the left hand term is purely imaginary and the right hand term is real.
\qed

\vspace{3mm} 
{\small
Author's address: \vspace{2mm}\\
Andrei Teleman, LATP, CMI,   Universit\'e de Provence,  39  Rue F.
Joliot-Curie, 13453 Marseille Cedex 13, France,  e-mail:
teleman@cmi.univ-mrs.fr. }

\end{document}